\documentclass{amsart}
\usepackage{amssymb, amsmath}
\usepackage{graphicx} 
\usepackage{multirow, array} 
\usepackage{float} 

\newtheorem{theorem}{Theorem}[section]

\newtheorem{remark}[theorem]{Remark}

\usepackage[dvipsnames]{xcolor}

\usepackage{pifont}

\newcommand{\spn}{\operatorname{span}}
\newcommand{\Ricci}{\operatorname{Ric}}
\newcommand{\bach}{\mathfrak{B}}

\newcommand{\tu}{\tilde u}
\newcommand{\te}{\tilde e}
\newcommand{\f}{f}
\newcommand{\tf}{\tilde f}

\begin{document}
\title[Conformally Einstein Lorentzian Lie groups]{Conformally Einstein Lorentzian Lie groups with Heisenberg symmetry}
\author[Calvi\~no-Louzao, Garc\'ia-R\'io, Guti\'errez-Rodr\'iguez, V\'azquez-Lorenzo]{E. Calvi\~no-Louzao, E. Garc\'ia-R\'io, I.  Guti\'errez-Rodr\'iguez, \\ R. V\'azquez-Lorenzo}
\address{ECL: I.E.S Eduardo Pondal, 15706 Santiago de Compostela,  Spain}
\email{estebcl@edu.xunta.gal}
\address{EGR: Department of Mathematics, CITMAga, University of Santiago de Compostela, 15782 Santiago de Compostela, Spain}
\email{eduardo.garcia.rio@usc.es}
\address{IDGR: Facultad de Ciencias de la Educación y del Deporte, University of Vigo, 36005 Pontevedra, Spain}
\email{ixchel.dzohara.gutierrez.rodriguez@uvigo.es}
\address{RVL: I.E.S de Ribadeo Dionisio Gamallo, 27700 Ribadeo,  Spain}
\email{ravazlor@edu.xunta.gal}
\thanks{Supported by projects PID2019-105138GB-C21 and PID2021-127075NA-I00 (AEI/FEDER, Spain) and ED431C 2019/10, ED431F 2020/04  (Xunta de Galicia, Spain).}
	
\begin{abstract}
We describe all Lorentzian semi-direct extensions of the Heisenberg group which are conformally Einstein. As a by side result,  
Bach-flat left-invariant Lorentzian metrics on semi-direct extensions of the Heisenberg group are classified, thus providing new background solutions in conformal gravity.
\end{abstract}
\maketitle

\section{Introduction}
Despite the great success of General Relativity, there are still
	some open issues regarding big distance scales which suggest the necessity of modifications in the theory of gravity itself.
	Conformal Geometry is one of the different approaches to modify gravity, which belongs to the category of higher-order derivatives in the action.
	
	The equation of motion for the metric in General Relativity is derived by a functional variation, $I=I_{EH}+I_\wedge+I_m$, of the Einstein-Hilbert action:
	$$
	I:g\mapsto I(g)=-\frac{1}{16\pi G}\int d^4x\sqrt{-g}(\tau-2\wedge)+I_m
	$$
	where $\tau$ is the scalar curvature, $I_m$ is the matter part and $I_\wedge$ denotes the action for the cosmological constant $\wedge$. The corresponding equations of motion are given by $(\rho-\frac{1}{2}\tau g)-\wedge g=-8\pi GT$, where  $T$ is the energy-momentum tensor, and $\rho$ is the Ricci tensor.
	It was originally expected that a modification of gravity should reproduce these equations, but this is not necessarily the case.
	
	Conformal gravity refers to gravity theories that are invariant under conformal transformations in the pseudo-Riemannian geometrical sense. The simplest theory in this category has the square of the norm of the Weyl tensor as the Lagrangian \cite{Mannheim1, Mannheim2}
	$$
	I:g\mapsto I(g)=\int d^4x\sqrt{-g}\,\| W\|^2+I_m,
	$$
	which has to be contrasted with the usual Einstein-Hilbert action, where the Lagrangian is just the scalar curvature.
	The equation of motion upon varying the metric is $\bach=0$, which is of fourth-order, where $\bach$ is the Bach tensor.

A four-dimensional pseudo-Riemannian manifold $(M^n,g)$ is said to be \emph{Bach-flat} if the Bach tensor
\begin{equation}\label{eq:EGR-1}
\bach=\operatorname{div}_2\operatorname{div}_4W+\tfrac{1}{2}W[\rho]
\end{equation}
vanishes identically, where $W$ denotes the Weyl conformal curvature tensor and $W[\rho]$ is the contraction of the Weyl tensor and the Ricci tensor $W[\rho]_{ij}=W_{iajb}\rho^{ab}$.
Locally conformally flat metrics have $W=0$ and thus they are trivially Bach-flat.
Einstein metrics are also Bach-flat since the Weyl tensor is traceless and $\operatorname{div}W=0$ in the Einstein case. 
The Bach tensor is the gradient of the quadratic curvature functional given by the $L^2$-norm of the Weyl tensor,  which is preserved by conformal transformations in dimension four \cite{Bach}. Hence four-dimensional conformally Einstein metrics are also Bach-flat.

A pseudo-Riemannian manifold $(M^n,g)$ is \emph{conformally Einstein} if there is an Einstein representative of the conformal class $[g]$.  Equivalently, a  conformally related metric $\overline{g}=\varphi^{-2}g$ is Einstein if and only if there exists a nowhere zero solution of the overdetermined PDE
\begin{equation}
\label{eq:hoy2}
(n-2)\operatorname{Hes}_\varphi+\varphi\rho=\tfrac{1}{n}\{(n-2)\Delta\varphi+\varphi\tau\} g ,
\end{equation}
where $\operatorname{Hes}_\varphi$ is the Hessian tensor of $\varphi$ and $\Delta\varphi=\operatorname{tr}_g\operatorname{Hes}_\varphi$ is the Laplacian of the (locally defined) conformal factor $\varphi$. Equation \eqref{eq:hoy2}, which is trivial in dimension two,  was originally considered by Brinkmann \cite{Bri1}. Despite its apparent simplicity, the integration of \eqref{eq:hoy2} is surprisingly difficult.

If a conformal metric $\overline{g}=e^{-2\sigma}g$ is Einstein, then $\overline{\operatorname{div}}\,\overline{W}=0$, and thus $\operatorname{div}_4W-W(\,\cdot\,,\,\cdot\,,\,\cdot\,,\nabla\sigma)=0$. 
A pseudo-Riemannian manifold is \emph{conformally Cotton-flat} if there is a conformal metric  $\overline{g}=e^{-2\sigma}g$ which is Cotton-flat. Equivalently, there is a function $\sigma$ so that 
\begin{equation}\label{eq:hoy3}
\operatorname{div}_4W-W(\,\cdot\,,\,\cdot\,,\,\cdot\,,\nabla\sigma)=0. 
\end{equation}
More generally one says that $(M,g)$ is a \emph{conformal $C$-space} if there is a (not necessarily gradient) vector field $X$ so that 
$\operatorname{div}_4W-W(\,\cdot\,,\,\cdot\,,\,\cdot\,,X)=0$ (see \cite{GN}).
The special significance of the conformally Cotton-flat property was given in \cite{KTN}, where it is shown that a weakly-generic Bach-flat manifold is conformally Einstein if and only if it is conformally Cotton-flat, where being weakly generic means that the Weyl tensor, viewed as a map $W:TM\rightarrow\otimes^3 TM$, is injective.

Spacetimes admitting  a null parallel vector field $\ell$ have been widely
studied in General Relativity, where they are called \emph{$pp$-waves} in the transversally 
flat case, i.e., if the curvature endomorphism satisfies $R(\ell^\perp,\ell^\perp)=0$ \cite{13}, where $R$ denotes the curvature tensor. Furthermore, the spacetime is a \emph{plane wave} if, in addition, the curvature tensor is transversally parallel (i.e., $\nabla_{\ell^\perp}R = 0$).  
Four-dimensional $pp$-waves were discovered in a mathematical context by Brinkmann \cite{Bri}. In Physics, plane waves and $pp$-waves appeared in General Relativity, where they play an important role (see \cite{BPR, 13}).
It was shown in \cite{BGV19} that four-dimensional $pp$-waves are conformally Einstein if and only if the Cotton tensor vanishes, or equivalently, $\operatorname{div}W=0$. (See also \cite{LN10} for a description of Bach-flat $pp$-waves). Consequently plane waves are conformally Einstein. 

Connected and simply connected four-dimensional Lie groups are either products $SU(2)\times\mathbb{R}$, $\widetilde{SL}(2,\mathbb{R})\times\mathbb{R}$, or one of the solvable semi-direct extensions of three-dimensional unimodular Lie groups $\widetilde E(2)\rtimes\mathbb{R}$, $E(1,1)\rtimes\mathbb{R}$, $H_3\rtimes\mathbb{R}$ or $\mathbb{R}^3\rtimes\mathbb{R}$, where $\widetilde E(2)$, $E(1,1)$, $H_3$ and $\mathbb{R}^3$ denote the simply connected Euclidean, Poincaré, Heisenberg and   Abelian three-dimensional Lie groups, respectively. 
Bach-flat left-invariant Riemannian metrics on four-dimensional Lie groups constitute a small class (see \cite{AGS, CL-GM-GR-GR-VL}). The situation is very different in the Lorentzian setting where there are plenty of Bach-flat left-invariant metrics. 
A general classification of these metrics seems quite an unfeasible task at this time. That is why in this paper we focus on a special family of Lie groups, the semi-direct extensions $H_3\rtimes\mathbb{R}$ of the Heisenberg group.  Hence we restrict to four-dimensional simply connected Lie groups which have the Heisenberg group as a normal subgroup acting with cohomogeneity one.
Left-invariant Einstein metrics on $H_3\rtimes\mathbb{R}$ were described in \cite{CaZa}, where it is shown that they are Ricci-flat plane waves, or of constant non-positive sectional curvature.
We analyze all left-invariant Lorentzian metrics and obtain classification results for the Bach-flat and conformally Einstein ones.
Since the Bach tensor and the conformally Einstein property are invariant by homotheties, we work at the homothetical level to simplify the discussion. It is important to emphasize that homotheties need not to be given by isomorphisms of the Lie groups. Hence in order to preserve the Lie group structure we work up to isomorphic homotheties.

\subsection{Summary of results}

Indecomposable Lorentzian symmetric spaces are irreducible (and hence of constant sectional curvature) or locally isometric to a Cahen-Wallach symmetric space \cite{Cahen}. Four-dimensional Cahen-Wallach symmetric spaces are isometric to $\mathbb{R}^4$ with coordinates $(u,v,x^1,x^2)$ and metric given by
$g=dudv+H(v,x^1,x^2)dvdv+dx^1dx^1+dx^2dx^2$, with $H(v,x^1,x^2)=\sum\lambda_i (x^i)^2$. Being plane waves,	Cahen-Wallach symmetric spaces are conformally Einstein and thus Bach-flat. Moreover, products of lower dimensional symmetric spaces are Bach-flat (indeed conformally Einstein) unless they are products of surfaces $N^2(c_1)\times N^2(c_2)$ with constant sectional curvature $c_1^2\neq c_2^2$.
Henceforth   we consider Lorentzian Lie groups $H_3 \rtimes \mathbb{R}$ which are not locally symmetric.

As usual, the Lie algebra $\mathfrak{h}_3$ of the Heisenberg group $H_3$ is described by a basis $\{ v_1,v_2,v_3\}$ with Lie bracket $[v_1,v_2]=v_3$. One-dimensional semi-direct extensions $\mathfrak{g}=\mathfrak{h}_3\rtimes\mathfrak{r}$ are determined by derivations of the Heisenberg algebra $\mathfrak{h}_3$, that becomes an ideal in $\mathfrak{h}_3\rtimes\mathbb{R}$. Furthermore if $\langle\cdot,\cdot\rangle$ is a Lorentzian inner product on $\mathfrak{g}$, its restriction to $\mathfrak{h}_3$ may be degenerate, of Riemannian signature, or of Lorentzian signature. This motivates a separate study of the three distinct, although not necessarily disjoint, situations above.

\begin{remark}\label{re:producto}
	\rm
	A special case of the analysis in this paper correspond to the \emph{oscillator algebra} and the corresponding \emph{oscillator group} (see, for example, \cite{MuRi,St}), and some particular metrics on semi-direct extensions $H_3\rtimes\mathbb{R}$ (see \cite{CaZa2}). 
	
	Another special situation is that of direct extensions $H_3\times\mathbb{R}$ so that the semi-direct extension reduces to the product Lie group.
	Rahmani showed in \cite{Rahmani} that there exist three non-homothetic classes of left-invariant Lorentzian metrics in the Heisenberg group $H_3$.  Kondo and Tamaru have recently shown in \cite{KT} that there exist exactly six non-homothetic classes of left-invariant Lorentzian metrics on $H_3\times\mathbb{R}$ up to automorphisms, which are described by 	the Lie algebra structures
	$$
	[e_1,e_2]=-(\alpha e_1-e_4),\quad
	[e_2,e_3]=\beta (\alpha e_1-e_4),\quad
	[e_2,e_4]=\alpha (\alpha e_1-e_4),
	$$ 
	where $\{ e_1,e_2,e_3, e_4\}$ is an orthonormal basis of $\mathfrak{h}_3\times\mathbb{R}$ with $e_4$ timelike, and the parameters $(\alpha,\beta)\in\{ (0,0), (1,0), (1,1), (2,0), (2,\sqrt{3}), (2,2)\}$.
	
	A direct calculation shows that $(\alpha,\beta)=(1,0)$ determines a flat metric. In the other cases the metric is Bach-flat if and only if $(\alpha,\beta)=(1,1)$, in which case it is locally conformally flat and locally symmetric, or  $(\alpha,\beta)=(2,\sqrt{3})$, in which case the Ricci operator and the Weyl curvature operator acting on the space of two-forms are two-step nilpotent. Furthermore, in the latter case, the Ricci tensor is parallel but the metric is not locally symmetric although the curvature tensor is modeled on a symmetric space.
	This metric corresponds to a plane wave whose derived algebra $[\mathfrak{g},\mathfrak{g}]$ is spacelike and the restriction of the metric to the center $\mathfrak{z}(\mathfrak{g})$ is of signature  $(+,-,0)=(1,0,1)$. 
	Bach-flat $pp$-wave left-invariant metrics on non-product semi-direct extensions $H_3\rtimes\mathbb{R}$ are discussed in $\S$\ref{sse:plane-wave bach flat}.
\end{remark}

\subsubsection{Non-trivial conformally Einstein metrics}\label{sse:conformally Einstein}
A conformally Einstein metric is said to be \emph{non-trivial} if it is neither Einstein, nor locally conformally flat, nor a plane wave.
The main result of this paper describes all  non-trivial conformally Einstein semi-direct Lorentzian extensions of the Heisenberg group as follows.

\begin{theorem}\label{th:ce}
	Let $(H_3 \rtimes \mathbb{R},\langle\cdot,\cdot\rangle)$ be a non-symmetric  semi-direct extension of the Heisenberg group equipped with a left-invariant Lorentzian metric. 
	Then, the metric  is non-trivial conformally Einstein if and only if it is isomorphically homothetic to one of the following:
	\begin{enumerate}
	\item[(D)] The restriction of the metric to $\mathfrak{h}_3$ is degenerate and the metric is determined by
	\begin{enumerate}
	\item[(D.i)] $[u_1,u_4]=   u_1$, 
	$[u_2,u_3]=  u_1$, 
	$[u_2,u_4]=   u_2$. 
	
	\item[(D.ii)] 
	$[u_1,u_3]=   u_1$, 
	$[u_1,u_4]= \varepsilon   u_1$,   
	$[u_2,u_3]= \alpha u_1$,  
	$[u_2,u_4]=   \varepsilon  u_2$,
	with $\varepsilon^2=1$  and $\alpha\geq 0$.
\end{enumerate}
Here  $\{u_1,u_2,u_3,u_4\}$ denotes a pseudo-orthonormal basis of $\mathfrak{h}_3\rtimes \mathfrak{r}$  with $\langle u_1,u_1\rangle=\langle u_2,u_2\rangle=\langle u_3,u_4\rangle=1$. 

\smallskip

\item[(R)] The restriction of the metric to $\mathfrak{h}_3$ is Riemannian and the metric is determined by
\[ 
[e_1,e_2]=  e_3,  \quad
[e_1,e_4]=\alpha e_1 +   e_3,   \quad 
[e_3,e_4]= \alpha  e_3, 
\quad\text{with}\quad
\alpha>0,
\]
where   $\{e_1,e_2,e_3,e_4\}$ denotes an orthonormal basis of $\mathfrak{h}_3\rtimes \mathfrak{r}$  with $e_4$ timelike.

\smallskip

\item[(L)] The restriction of the metric to $\mathfrak{h}_3$ is Lorentzian and the metric is determined by
\begin{enumerate}
\item[(L.i)]
	$[e_1,e_3]=-\alpha e_2$,
	$[e_1,e_4]=e_1 + \alpha e_2$,
	$[e_2,e_4]=  e_2$,
	where $\alpha>0$ and $\{e_1,e_2,e_3,e_4\}$ denotes an  orthonormal basis of $\mathfrak{h}_3\rtimes \mathfrak{r}$ with $e_3$ timelike.

\item[(L.ii)]
	$[u_1,u_3]=-u_2$,
	$[u_1,u_4]=4 u_1+\alpha u_2$,
	$[u_2,u_4]=4 u_2$,
	where $\alpha\in\mathbb{R}$ and  
$\{u_1,u_2,u_3,u_4\}$  denotes a  pseudo-orthonormal basis of $\mathfrak{h}_3\rtimes \mathfrak{r}$ with $\langle u_1,u_2\rangle=\langle u_3,u_3\rangle=\langle u_4,u_4\rangle=1$.	
\item[(L.iii)]
	$[u_1,u_3]=- u_2$,
	$[u_1,u_4]=-u_1$,
	$[u_2,u_4]=3  u_2$, 
	$[u_3,u_4]=  4  u_3,
$
where
$\{u_1,u_2,u_3,u_4\}$  denotes a  pseudo-orthonormal basis of $\mathfrak{h}_3\rtimes \mathfrak{r}$, with $\langle u_1,u_2\rangle=\langle u_3,u_3\rangle=\langle u_4,u_4\rangle=1$.
\end{enumerate}
\end{enumerate}
Moreover, in all cases but (L.iii) the left-invariant metrics are conformally equivalent to a Ricci-flat $pp$-wave.
\end{theorem}

Considering the eigenvalue structure of the Ricci operator for the different cases in Theorem~\ref{th:ce} given in Remark~\ref{re:degenerate case}, Remark~\ref{re:Riemann case}, and Remark~\ref{re:Lorentz case1}, one has that all of them correspond to different homothetical classes except possibly (D.ii) and (R) for $\varepsilon=-1$ (resp., (D.ii) and (L.i) for $\varepsilon=1$). Moreover, a direct calculation of $W[\rho]$ and the eigenvalues of the associated $(1,1)$-tensor field shows that the above mentioned cases cannot be homothetic.
Alternatively, the non-existence of homotheties between the classes above also follows from the work of \cite{CoHePe}, just considering the corresponding orthonormal bases of eigenvectors of the Ricci operators.

\begin{remark}\label{re:degenerate case}\rm
	The Ricci operator of the left-invariant metric in Theorem~\ref{th:ce}-(D.i) is diagonalizable with eigenvalues $\{1$, $1$, $-1$, $-1\}$, so that the scalar curvature $\tau=0$.
	The Weyl curvature operator $W:\Lambda^2\rightarrow\Lambda^2$ is two-step nilpotent  and hence not   weakly generic.
	
	The space of scalar quadratic curvature invariants of a pseudo-Riemannian manifold is generated by $\{\tau^2,\|\rho\|^2,\| R\|^2,\Delta\tau\}$. Hence any quadratic curvature functional is a linear combination of the corresponding $L^2$-norms: $g\mapsto \Phi_{abc}(g)=\int dx^4\sqrt{-g}\,\,\{ a\|\rho\|^2+b\| R\|^2+c\tau^2\}$. Due to the four-dimensional Gauss-Bonnet Theorem, any quadratic curvature functional in dimension four is equivalent to one of 
	$$
	\mathcal{S}:g\mapsto\mathcal{S}(g)=\int dx^4\sqrt{-g}\,\tau^2 
	\quad\text{or}\quad
	\mathcal{F}_t:g\mapsto\mathcal{F}_t(g)=\int dx^4\sqrt{-g}\,\left\{\|\rho\|^2+t\tau^2\right\}\,.
	$$
	In particular, one has that the functional given by the $L^2$-norm of the Weyl conformal curvature tensor is equivalent to $\mathcal{F}_{-1/3}$.
	Einstein metrics are critical for all quadratic curvature functionals in dimension three and four (but not necessarily in higher dimensions). Moreover, a four-dimensional metric is critical for all quadratic curvature functionals if and only if it is critical for two-distinct quadratic curvature functionals. In particular, any Bach-flat metric with vanishing scalar curvature is critical for all quadratic curvature functionals.
	Since plane waves are Bach-flat and have two-step nilpotent Ricci operator, they are critical for all quadratic curvature functionals.
	The metric in Theorem~\ref{th:ce}-(D.i) has vanishing scalar curvature.  Hence it is critical for all quadratic curvature functionals although it is not a plane wave nor locally symmetric (even not modeled on a symmetric space). 
	
	\medskip
	
	The Ricci operator of metrics in Theorem~\ref{th:ce}-(D.ii) has four different real eigenvalues $\{ \varepsilon(-1\pm\sqrt{\alpha^2+2})$,   $\varepsilon(-2\pm\sqrt{\alpha^2+1}) \}$ so that the scalar curvature $\tau=-6\varepsilon$. Moreover, fixing $\varepsilon$, $\alpha\geq 0$ determines the homothetic class, since the homothetic invariant $\tau^{-2}\|\rho\|^2 = \frac{1}{9}(\alpha^2+4)$.
	Metrics are not weakly generic since the Weyl curvature operator acting on the space of two-forms is two-step nilpotent.
\end{remark}

\begin{remark}\label{re:Riemann case}\rm
	The Ricci operator  of left-invariant metrics in Theorem~\ref{th:ce}-(R) has four different real eigenvalues $\{ \alpha (2\alpha\pm 1)$, $\alpha(\alpha\pm\sqrt{\alpha^2+1})\}$, so that the scalar curvature $\tau=6\alpha^2$.
	Moreover,  $\alpha> 0$ determines different homothetic classes since the homothetic invariant $\tau^{-2}\|\rho\|^2 = \frac{1}{9}(\frac{1}{\alpha^2}+3)$.
	Metrics are	not weakly generic since the Weyl curvature operator acting on the space of two-forms is two-step nilpotent.
\end{remark}

\begin{remark}\label{re:Lorentz case1}\rm
	The Ricci operator of left-invariant metrics in Theorem~\ref{th:ce}-(L.i) has four different real eigenvalues $\{-2\pm\alpha$, $ -1\pm\sqrt{\alpha^2+1} \}$, so that the scalar curvature $\tau=-6$. The metric is not   weakly generic since the Weyl curvature operator is two-step nilpotent.
	Moreover, the parameter  $\alpha> 0$ determines the homothetic classes since $\tau^{-2}\|\rho\|^2 = \frac{1}{9}(\alpha^2+3)$. 

	Metrics in Theorem~\ref{th:ce}-(L.ii) have Ricci operator with eigenvalues $\{0$, $-32$, $-32$, $-32\}$, so that the scalar curvature $\tau=-96$. The Ricci operator is diagonalizable if $\alpha=0$, having  a double root of the minimal polynomial otherwise. Moreover, the metrics are not weakly generic since the Weyl curvature operator acting on the space of two-forms is two-step nilpotent.
		
	The Ricci operator of the left-invariant metric in Theorem~\ref{th:ce}-(L.iii) is diagonalizable,  $\Ricci=\operatorname{diag}[-6$, $-6$, $-24$, $-18]$ and therefore $\tau=-54$.
	The Weyl curvature operator has eigenvalues $\{-4,2,2,-4,2,2\}$, where the eigenvalue $2$ is a double root of the minimal polynomial. Hence, it is weakly generic.
\end{remark}

\begin{remark}
	\rm\label{re:solitonesalgebraicos}
	The two-loop renormalization group flow (RG2 flow for short)  is a perturbation of the Ricci flow $\partial_tg_t=-2\rho[g_t]$, which mathematically is described by
	$\partial_tg_t=-2RG[g_t]$. The symmetric $(0,2)$-tensor field $RG=\rho +\frac{\Upsilon}{4}\check{R}$, where $\Upsilon$ denotes a positive coupling constant and  $\check{R}$ is the symmetric $(0,2)$-tensor field given by $\check{R}_{ij}=R_{iabc}R_j{}^{abc}$.
	We refer to \cite{Carfora, GGI2} and references therein for more information on the RG2 flow.
	Genuine fixed points of the flow are provided by those manifolds where the   tensor field $RG$ vanishes, i.e., $\rho +\frac{\Upsilon}{4}\check{R}=0$.
	Given a one-parameter family $\psi_t$ of diffeomorphisms of $M$ (with $\psi_0=\operatorname{Id}$), a solution of the form $g(t)=\sigma(t)\psi_t^*g$ (where $\sigma$ is a real-valued function with $\sigma(0)=1$) is said to be a \emph{self-similar solution}.
	A triple $(M,g,X)$, where $X$ is a vector field on $M$, is called an  RG2 \emph{soliton} if $\mathcal{L}_Xg+RG=\lambda g$ for some $\lambda\in\mathbb{R}$. Further the soliton is said to be \emph{expanding}, \emph{steady} or \emph{shrinking} if $\lambda <0$, $\lambda=0$, or $\lambda >0$, respectively. 
	Any self-similar solution of the RG2 flow is an  RG2 soliton just considering the vector field $X$ generated by the one-parameter group of diffeomorphisms $\psi_t$.
	Since the two terms comprising $RG$ behave differently under homotheties ($\rho[\kappa g]=\rho[g]$ and $\check{R}[\kappa g]=\frac{1}{\kappa}\check{R}[g]$), one has that the converse holds only for steady solitons, in which case $\psi_t$ is the one-parameter group of diffeomorphisms associated to the vector field $X$ determined by the soliton equation $\mathcal{L}_Xg+RG=0$ and $g(t)=\psi_t^*g$ is a self-similar solution (see~\cite{Wears}). 
	
	Let $G$ be a Lie group with left-invariant metric $\langle\cdot,\cdot\rangle$ and let $(\mathfrak{g},\langle\cdot,\cdot\rangle)$ denote the corresponding Lie algebra.  An  RG2 \emph{algebraic soliton} is a derivation of the Lie algebra $\mathfrak{g}$ given by $\mathfrak{D}=\widehat{RG}-\lambda\operatorname{Id}$, where $\widehat{RG}$ is the $(1,1)$-tensor field metrically equivalent to $RG$ and $\lambda\in\mathbb{R}$. RG2 algebraic solitons give rise to RG2 solitons (where the vector field $X$ is associated to a one-parameter group of automorphisms of $G$ determined by the derivation $\mathfrak{D}$) as in the Ricci flow case (see~\cite{Lauret, Wears}).
	
	Now, a straightforward calculation shows that left-invariant metrics in Theorem~\ref{th:ce}-(L.ii) are steady  RG2 algebraic solitons with  $RG2=\rho+\frac{1}{32}\check{R}$. Therefore they are steady RG2 solitons and thus also self-similar solutions of the flow with two-step nilpotent tensor field $\widehat{RG}$. Moreover, they are not algebraic Ricci solitons.
	
	In contrast with the previous situation, the left-invariant metric in Theorem~\ref{th:ce}-(L.iii) is a shrinking RG2 algebraic soliton  with  $RG2=\rho+\frac{1}{10}\check{R}$.  Therefore they are shrinking RG2 solitons with diagonal tensor field $\widehat{RG}$, but not algebraic Ricci solitons.
\end{remark}

\subsubsection{Strictly Bach-flat metrics}\label{sse:Bach flat}

We say that a Bach-flat metric is \emph{strict} if it is not conformally Einstein (and thus not a  plane wave).
In contrast with Theorem~\ref{th:ce} and the case of $pp$-waves in Theorem~\ref{th:plane-waves}, strictly Bach-flat semi-direct extensions $H_3\rtimes\mathbb{R}$ are quite rare and they are described in the following 

\begin{theorem}\label{th:strictly Bach flat} 
	Let $(H_3 \rtimes \mathbb{R},\langle\cdot,\cdot\rangle)$ be a non-symmetric  semi-direct extension of the Heisenberg group equipped with a left-invariant Lorentzian metric.  
	Then, the metric  is  strictly Bach-flat if and only if it is isomorphically homothetic to  one of the following Lie algebras:
	\begin{enumerate}
		\item[(i)] 
		$[v_1,v_2] = -3 \sqrt{14} v_2 + 2 \varepsilon  \sqrt{11} v_3$,
		$[v_1,v_3] = \sqrt{14}v_3$,
		$[v_1,v_4] = -2\sqrt{14} v_4$, 
		$[v_2,v_3] = 4 v_4$,
		where $\varepsilon^2=1$ and $\{v_1,v_2,v_3,v_4\}$ denotes an orthonormal basis of $\mathfrak{h}_3\rtimes \mathfrak{r}$ with $v_1$ and $v_4$ spacelike and 
		$\langle v_2,v_2\rangle=-\langle v_3,v_3\rangle=\bar\varepsilon=\pm1$. 
		
		\smallskip
		
		\item[(ii)]
		$[u_1,u_3] = 			
		-\tfrac{1}{\sqrt{-(\alpha+\beta^2)}}u_2$,
		$[u_1,u_4] = 			
		\tfrac{1}{\sqrt{-(\alpha+\beta^2)}}
		(\beta u_1+\gamma u_2+ u_3)$,
		
		\smallskip
		\noindent
		$[u_3,u_4] = 			
		\tfrac{1}{\sqrt{-(\alpha+\beta^2)}}
		\left(
		\alpha u_1
		-\frac{\alpha+3\beta^2}{4\alpha} u_2
		-\beta u_3
		\right)$,
		where  the constants $\alpha$, $\beta$, $\gamma$ satisfy
		$\gamma =\tfrac{1}{4\alpha^2} \left(\!
		(3\alpha +\beta^2)\beta
		\pm 2  \sqrt{-4\alpha^4 -2 (\alpha+\beta^2  )^3} 
		\right)$,
		with  
		$-\frac{1}{2}\!\leq\!\alpha\!<\!0$ and
		$\beta^2\leq -\alpha\left(
		\sqrt[3]{2\alpha}+1\right)$. Moreover
	$\{u_1,u_2,u_3,u_4\}$  denotes a  pseudo-orthonormal basis of $\mathfrak{h}_3\rtimes \mathfrak{r}$ with $\langle u_1,u_2\rangle=\langle u_3,u_3\rangle=\langle u_4,u_4\rangle=1$.	
		\end{enumerate}
\end{theorem}

\begin{remark}\label{re:Lorentz case2}\rm
	The Ricci operator of left-invariant metrics in Theorem~\ref{th:strictly Bach flat}-(i) is diagonalizable, $\operatorname{Ric}= \operatorname{diag}[  -174$, $-138$, $42$, $-120]$, so $\tau=-390$. Moreover, $\bar\varepsilon$ determines the spacelike or timelike character of the eigenvectors associated to the eigenvalues $-138$ and $42$.
	Hence the parameter $\bar\varepsilon$ determines two non-homothetic classes. Furthermore the Weyl curvature operator $W:\Lambda^2\rightarrow\Lambda^2$ has six-distinct non-zero complex eigenvalues and thus it is weakly generic.

	The characteristic polynomial of the Ricci operator corresponding to metrics in Theorem~\ref{th:strictly Bach flat}-(ii) is given by
	\[
	\operatorname{det}(\operatorname{Ric}-\lambda\operatorname{Id})=
	\lambda^4 
	-\tfrac{3}{4} \lambda^3
	+ \tfrac{3}{16} \lambda^2
	-\tfrac{9}{64} \lambda
	+ \tfrac{\alpha^4}{ 8(\alpha + \beta^2)^3}.
	\]
A straightforward calculation shows that the discriminant, given by
	\[
	\begin{array}{l}
	\Delta=\tfrac{1}{262144 (\alpha+\beta^2)^9}\left(
	131072 \alpha^{12}
	-74160   \alpha^8(\alpha+\beta^2)^3\right.
	\\
	\noalign{\medskip}
	\phantom{\Delta=\tfrac{1}{262144 (\alpha+\beta^2)^9}\left(\right.}
	\left.
	+  12312  \alpha^4 (\alpha+\beta^2)^6
	-2187 (\alpha+\beta^2)^9
	\right) ,
	\end{array}
	\]
	is  strictly negative considering the restrictions on $\alpha,\beta$ in Theorem~\ref{th:strictly Bach flat}. As a consequence,    the Ricci operator has   two real roots (with opposite sign) and two complex conjugate roots.
Moreover, the scalar curvature is strictly positive,  $\tau=\frac{3}{4}$. The Weyl curvature operator acting on the space of two-forms is three-step nilpotent. 
\end{remark}

\subsubsection{Bach-flat $pp$-wave metrics}\label{sse:plane-wave bach flat}
It was shown in Remark~\ref{re:producto} that all Bach-flat left-invariant metrics on the product $H_3\times\mathbb{R}$ are plane waves. Bach-flat left-invariant metrics on non-product semi-direct extensions $H_3\rtimes\mathbb{R}$ which are $pp$-waves are now given as follows.
It turns out that they are plane waves. We refer to \cite{BO03} for a classification of homogeneous plane waves.
	
\begin{theorem}\label{th:plane-waves}
	Let $(H_3 \rtimes \mathbb{R},\langle\cdot,\cdot\rangle)$ be a non-product semi-direct extension of the Heisenberg group equipped with a left-invariant Lorentzian metric.  
	 	If $(H_3 \rtimes \mathbb{R},\langle\cdot,\cdot\rangle)$  is a 	
		non-symmetric  Bach-flat $pp$-wave which is not locally conformally flat, then it is isomorphically homothetic to  one of the following:
		\begin{enumerate}
			\item[(D.i)] 
			The left-invariant metric determined by 
			\[
			\begin{array}{ll}
			[u_1,u_2]=  u_3 , &
			[u_1,u_4]= \kappa_1 u_1 + \kappa_2 u_2 +\kappa_4 u_3,
			\\
			\noalign{\medskip}
			[u_2,u_4]=  -\kappa_2 u_1 + \kappa_3 u_2+\kappa_5 u_3, &
			[u_3,u_4]=(\kappa_1+\kappa_3) u_3 ,
			\end{array}
			\]
			where  $\kappa_1$, $\kappa_2$, $\kappa_3$, $\kappa_4$, $\kappa_5\in\mathbb{R}$, with $(\kappa_1-\kappa_3)(2\kappa_2+1)\neq 0$.
			
			\smallskip
			
			\item[(D.ii)]
			The left-invariant metric determined by
			$$
			[u_1,u_3]=u_1,\quad
			[u_2,u_3]=  \kappa\,  u_1,\quad \kappa> 0,
			$$
		\end{enumerate}
		where  $\{u_1,u_2,u_3,u_4\}$ is a pseudo-orthonormal basis of $\mathfrak{h}_3\rtimes \mathfrak{r}$, with $\langle u_1,u_1\rangle=\langle u_2,u_2\rangle=\langle u_3,u_4\rangle=1$.
		
		\smallskip
		
		\begin{enumerate}	
			\item[(L)]
			The left-invariant metric determined by 
			\[ 
			[u_1,u_3]=-\varepsilon u_2,\quad
			[u_1,u_4]= \kappa_1 u_2+\kappa_2 u_3,\quad
			[u_3,u_4]= \kappa_3 u_2,
			\]
		\end{enumerate}	
		where  $\{u_1,u_2,u_3,u_4\}$ is a pseudo-orthonormal basis of $\mathfrak{h}_3\rtimes \mathfrak{r}$, with $\langle u_1,u_2\rangle=\langle u_3,u_3\rangle=\langle u_4,u_4\rangle=1$, and $\kappa_1,\kappa_2,\kappa_3\in\mathbb{R}$ with $\kappa_2(\kappa_2+\kappa_3)\neq0$, and $\varepsilon=\pm1$.
		
		\smallskip
		
		Moreover, all the cases above are conformally Einstein plane waves.
	\end{theorem}
	
	\begin{remark}\label{re:casos no producto}
		\rm
		Recall from Remark~\ref{re:producto} that all Bach-flat metrics on $H_3\times\mathbb{R}$ have parallel Ricci tensor. A straightforward calculation shows that the Ricci tensor of metrics corresponding to (D.i) is parallel if and only if $(\kappa_1+\kappa_3)(4\kappa_1\kappa_3+1)=0$. Hence in the generic situation these metrics are not homothetic to any left-invariant metric on the product $H_3\times\mathbb{R}$.

		The Lie group in (D.ii) is not isomorphic to the product $H_3\times\mathbb{R}$, although the Ricci tensor is parallel. Indeed, the underlying Lie group is solvable but not two-step nilpotent in contrast with  the product $H_3\times\mathbb{R}$.
		
		The Lie group in (L) is three-step nilpotent. Hence left-invariant metrics in (L) are  not isomorphically homothetic  to any left-invariant metric on the product $H_3\times\mathbb{R}$, although the Ricci tensor is parallel.
		
		A straightforward calculation shows that the sectional curvature of metrics corresponding to (D.i) does not depend on the parameters $\kappa_4$ and $\kappa_5$. Hence it follows from the work in \cite{Kulkarni} that a left-invariant metric (D.i) is homothetic (although not necessarily isomorphically homothetic) to a metric with $\kappa_4=\kappa_5=0$.
		Analogously, the sectional curvature of metrics (L) does not depend both on $\kappa_1$ and $\varepsilon$. Hence one has that a metrics in (L) is homothetic (but not isomorphically homothetic) to a metric with $\varepsilon=-1$ and $\kappa_1=0$.
	\end{remark}

\subsection{Bach tensor of left-invariant metrics and Gr\"obner bases}\label{se:Grobner}
Let $(G,\langle\cdot,\cdot\rangle)$ be a four-dimensional Lorentzian Lie group. It is now immediate that the Bach-flatness condition $\bach=0$ equals to a system of polynomial equations on the structure constants (given by the components $\bach_{ij}$ of the Bach tensor) which one has to solve in order to obtain a complete classification. When the system under consideration is simple, it is an elementary problem to find all common roots, but if the number of equations, unknowns and their degrees increase, it may become a quite unmanageable task. 

Given a set $\mathcal{S}$ of polynomials $\mathfrak{B}_{ij} \in \mathbb{R}[x_1, \dots, x_n]$, an $n$-tuple of real numbers $\vec{a}=(a_1, \dots, a_n)$ is a solution of $\mathcal{S}$ if and only if $\mathfrak{B}_{ij}(\vec{a}) = 0$ for all $i$, $j$. It is a fundamental observation to recognize that $\vec{a}$ is a solution of $\mathcal{S}$ if and only if it is a solution of $\mathcal{I} = \langle \mathfrak{B}_{ij} \rangle$, the ideal generated by the $\mathfrak{B}_{ij}$'s: if two sets of polynomials generate the same ideal, the corresponding zero sets must be identical. 
The theory of Gröbner bases provides a well-known strategy to solve rather large polynomial systems obtaining ``better'' polynomials that belong to the ideal generated by the initial polynomial system (see \cite{Cox} for more information on Gröbner bases).

Let $x^{\alpha}=x_1^{\alpha_1}\cdots x_n^{\alpha_n}$ with $\alpha\in\mathbb{Z}^n_{\geq 0}$ be a monomial in $\mathbb{R}[x_1, \dots, x_n]$. 
	A monomial ordering is any relation on the set of monomials $x^\alpha$ with $\alpha\in\mathbb{Z}^n_{\geq 0}$ satisfying
	\begin{enumerate}
		\item It is a total ordering on $\mathbb{Z}^n_{\geq 0}$.
		\item If $\alpha>\beta$ and $\gamma\in\mathbb{Z}^n_{\geq 0}$, then $\alpha+\gamma>\beta+\gamma$.
		\item $\mathbb{Z}^n_{\geq 0}$ is well-ordered, so that every non-empty subset of $\mathbb{Z}^n_{\geq 0}$ has a smallest element with respect to the given ordering.
	\end{enumerate}
	Establishing an ordering on $\mathbb{Z}^n_{\geq 0}$ will induce an ordering on the monomials. 
	For our purposes we will   use the \emph{lexicographical order} and the \emph{graded reverse lexicographical order}.  We say that $\alpha>_{lex}\beta$ if in the vector $\alpha-\beta\in\mathbb{Z}^n$  the leftmost non-zero entry is positive
	and we say that $\alpha>_{grevlex}\beta$ if $|\alpha|>|\beta|$ or $|\alpha|=|\beta|$     
	and the rightmost non-zero entry of $\alpha-\beta\in\mathbb{Z}^n$  is negative.
	We would like to emphasize that the Gröbner basis construction is very sensitive to the ordering on the variables and the monomials. For a certain ordering, a simple Gröbner basis can be obtained with a reduced number of polynomials, while for other orderings both the number of polynomials and their form can be completely unmanageable. 
	Lexicographical order for the monomials is the most appropriate in most cases to get simple bases. However, it is not always possible to use such ordering by computational reasons and other orderings must be taken into consideration.
	We therefore emphasize in each case the ordering under consideration for the monomials and the variables.
	
	Finally, it is worth to note that Gröbner bases are not unique for a given ordering on the variables and the monomials, since they may also depend on the algorithm used in the calculations. However, they are suitable to check the ideal membership problem and to decide whether a given polynomial is in the ideal under consideration. 
All the calculations in this paper have been done with {\sc Singular} \cite{Singular} and doubly checked with {\sc Mathematica}. The corresponding files are available from the Authors under request.

\subsection{Schedule of the paper}
We analyze the existence of Bach-flat metrics in semi-direct extensions $H_3\rtimes\mathbb{R}$ in Sections $\S$\ref{se:2}, $\S$\ref{se:Riemann} and $\S$\ref{se:EGR-Lorentzian}, depending on whether the induced metric in $\mathfrak{h}_3$ is degenerate, Riemannian or Lorentzian, respectively. The proof of Theorem~\ref{th:plane-waves} follows directly from this analysis and Remark~\ref{re:producto}. In all the cases we pay special attention to obtain simpler descriptions of the corresponding Bach-flat metrics by using suitable isomorphisms within the homothetic class. As a consequence the conformally Einstein equation \eqref{eq:hoy2} coupled with the conformally Cotton-flat equation \eqref{eq:hoy3} become tractable, and Theorem~\ref{th:ce} and Theorem~\ref{th:strictly Bach flat} are finally proven in Section~\ref{se:proof-1} and Section~\ref{se:proof-2}, respectively.

\section{Semi-direct extensions with degenerate normal subgroup $H_3$}\label{se:2}

In this section we analyze left-invariant Lorentzian metrics which are extensions of the three-dimensional unimodular Lie group $H_3$ equipped with a degenerate metric. 
Hence, let $\mathfrak{g}=\mathfrak{h}_3\rtimes\mathfrak{r}$ be a four-dimensional Lie algebra with a Lorentzian inner product $\langle\cdot,\cdot\rangle$  which restricts to a degenerate inner product on the subalgebra $\mathfrak{h}_3$. Let $\mathfrak{h}_3'=\spn\{ v\}$  be the derived subalgebra of $\mathfrak{h}_3$, $\mathfrak{h}_3'= [\mathfrak{h}_3,\mathfrak{h}_3]$. Since  the restriction of the metric to $\mathfrak{h}_3$ has signature  $(+,+,0)$ the vector $v$ may be spacelike or null (see \cite{CC}).   
Next we  analyze the vanishing of the Bach tensor in  those two cases by separate.

\subsection{\bf $\boldsymbol{\mathfrak{h}'_3=\operatorname{\bf span}\{v\}}$ is a null subspace}\label{se:deg-null subspace}

Setting $u_3=v$ we can take a pseudo-orthonormal basis $\{u_1,u_2,u_3,u_4\}$ of $\mathfrak{g}=\mathfrak{h}_3\rtimes \mathfrak{r}$, with $\langle u_1,u_1\rangle=\langle u_2,u_2\rangle=\langle u_3,u_4\rangle=1$,  so that  $\mathfrak{h}_3=\spn\{u_1,u_2,u_3\}$ and $\mathfrak{r}=\spn\{u_4\}$. Since $\operatorname{ad}(u_4)$ acts on $\mathfrak{h}_3$ as a derivation, one has the Lie brackets 
\[
[u_1,u_2]=\lambda_1 u_3, \quad
[u_1,u_3]=\lambda_2 u_3,  \quad
[u_2,u_3]=\lambda_3 u_3,  \quad
\underset{\underset{(i=1,2,3)}{}
}{[u_i,u_4]}=\sum_{j=1}^3 \alpha_i^j u_j ,
\]
for certain $\alpha_i^j\in\mathbb{R}$, where at least one of $\lambda_1$, $\lambda_2$ and $\lambda_3$ is non-zero. Next,  considering the endomorphism determined by the matrix $(\alpha_i^j)$, or equivalently the Jacobi identity, one is led to the following different possibilities   depending on the $\lambda_i$'s.

We will show that Bach-flat left-invariant metrics in this case are locally conformally flat or plane waves, thus not providing new examples of strictly Bach-flat structures.

\smallskip

\subsubsection{Case $\lambda_2=\lambda_3=0$}\label{se:deg-1-1}
If $\lambda_2=\lambda_3=0$, then necessarily $\lambda_1\neq 0$ and  
\[
\begin{array}{ll}
	[u_1,u_2]= \lambda_1 u_3 , &
	[u_1,u_4]= \gamma_1 u_1 + \gamma_2 u_2 + \gamma_3 u_3 ,
	\\
	\noalign{\medskip}
	[u_2,u_4]=  \gamma_4 u_1 + \gamma_5 u_2 + \gamma_6 u_3 , &
	[u_3,u_4]=(\gamma_1+\gamma_5) u_3 ,
\end{array}
\]
where $\gamma_1$, $\dots$, $\gamma_6\in\mathbb{R}$.
In this case, without further assumptions,  a direct  calculation shows that the Bach tensor vanishes. 
Moreover, $u_3$ is a null recurrent vector field and  the curvature tensor satisfies $R(x,y)=0$ and $\nabla_xR=0$ for all $x,y\in u_3^\perp=\spn\{u_1,u_2,u_3\}$. Finally,  the Ricci tensor of the  above metric  is determined by  $\rho_{44}=\frac{1}{2}\{ \lambda_1^2 + 4\gamma_1\gamma_5 -(\gamma_2+\gamma_4)^2\}$, which implies that the Ricci operator is isotropic  and therefore the underlying structure is a plane wave (see \cite{Leistner}).

Rotating the spacelike vectors $\{ u_1,u_2\}$ one may assume that $\gamma_4=-\gamma_2$.
Furthermore, rescaling the vectors $u_k$ by $\frac{1}{\lambda_1}u_k$ one may set $\lambda_1=1$ remaining in the same homothety class.  
Finally a straightforward calculation shows that the metric is locally conformally flat if and only if 
$(\gamma_1-\gamma_5)(2\gamma_2+1)= 0$, and the metric is locally symmetric if and only if 
 $\gamma_1=\gamma_5=0$, or $2\gamma_2+1=0$ and $(\gamma_1+\gamma_5)(4\gamma_1\gamma_5+1)=0$. 
It corresponds to Assertion (D.i) in Theorem~\ref{th:plane-waves}.

\smallskip

\subsubsection{Case $\lambda_2=0$,  $\lambda_3\neq 0$}\label{se:deg-1-2}

In this case, the Lie algebra structure is given by 
\[
\begin{array}{lll}
	[u_1,u_2]= \lambda_1 u_3 , &
	[u_1,u_4]= \gamma_1\lambda_3 u_1 + (\gamma_1-\gamma_2)\lambda_1 u_3  , &
	[u_2,u_3]= \lambda_3 u_3 ,
	\\
	\noalign{\medskip}
	[u_2,u_4]=  \gamma_3 u_1 + \gamma_4 u_3 , &
	[u_3,u_4]= \gamma_2\lambda_3  u_3 ,
\end{array}
\]  
where $\gamma_1$, $\dots$, $\gamma_4\in\mathbb{R}$. A straightforward calculation shows that the only non-zero component of the Bach tensor, $\mathfrak{B}_{ij}=\mathfrak{B}(u_i,u_j)$, corresponds to $\bach_{44}=\frac{1}{4}(\gamma_1^2\lambda_3^2+\gamma_3^2)\lambda_3^2$. Hence, since $\lambda_3\neq 0$, necessarily $\gamma_1=\gamma_3=0$. Finally, one directly checks  that the metric is locally conformally flat and locally symmetric.

\smallskip

\subsubsection{Case $\lambda_2\neq 0$}\label{se:deg-1-3}

If $\lambda_2\neq 0$,  then one has
\[
\begin{array}{l}
	[u_1,u_2]= \lambda_1 u_3 , \quad
	[u_1,u_3] = \lambda_2 u_3, \quad
	[u_2,u_3]= \lambda_3 u_3, \quad
	[u_3,u_4]= \gamma_4\lambda_2  u_3,
	\\
	\noalign{\medskip}
	[u_1,u_4]= -\gamma_1\lambda_2\lambda_3 u_1 +\gamma_1\lambda_2^2 u_2 +\gamma_2\lambda_2 u_3 , 
	\\
	\noalign{\medskip}
	[u_2,u_4]=  -\gamma_3\lambda_3 u_1 +\gamma_3\lambda_2 u_2 
	+( \gamma_1\lambda_1\lambda_3-(\gamma_3-\gamma_4)\lambda_1+\gamma_2\lambda_3  ) u_3  ,
\end{array}
\]  
where $\gamma_1$, $\dots$, $\gamma_4\in\mathbb{R}$. As in the previous case, the Bach tensor is determined by just one component, $\bach_{44}=\frac{1}{4}(\gamma_1^2\lambda_2^2+\gamma_3^2)(\lambda_2^2+\lambda_3^2)^2$, which implies $\gamma_1=\gamma_3=0$. Again, a direct calculation shows that the metric is locally conformally flat and locally symmetric.

\subsection{\bf $\boldsymbol{\mathfrak{h}'_3=\spn\{v\}}$ is a spacelike subspace}\label{se:deg-spacelike subspace}

In this case, we set  $u_1=\frac{v}{\|v\|}$ and consider a pseudo-orthonormal basis $\{u_1,u_2,u_3,u_4\}$ of $\mathfrak{g}=\mathfrak{h}_3\rtimes \mathfrak{r}$, with $\langle u_1,u_1\rangle=\langle u_2,u_2\rangle=\langle u_3,u_4\rangle=1$,  where  $\mathfrak{h}_3=\spn\{u_1,u_2,u_3\}$ and $\mathfrak{r}=\spn\{u_4\}$, so that 
\[
[u_1,u_2]=\lambda_1 u_1, \quad
[u_1,u_3]=\lambda_2 u_1,  \quad
[u_2,u_3]=\lambda_3 u_1,  \quad
\underset{\underset{(i=1,2,3)}{}
}{[u_i,u_4]}=\sum_{j=1}^3 \alpha_i^j u_j ,
\]
for certain $\alpha_i^j\in\mathbb{R}$ and where at least one of $\lambda_1$, $\lambda_2$ and $\lambda_3$ is non-zero.
We proceed as in Section~\ref{se:deg-null subspace}, to have the following different possibilities depending on the $\lambda_i$'s.

\smallskip

\subsubsection{Case $\lambda_1=\lambda_2=0$}\label{se:deg-2-1}

In this case, necessarily  $\lambda_3\neq 0$ and     the Lie algebra structure is given by
\[
\begin{array}{l}
	[u_1,u_4]= \gamma_1 u_1, \quad
	[u_2,u_3]= \lambda_3 u_1, \quad
	[u_2,u_4]= \gamma_2 u_1 +\gamma_3 u_2 + \gamma_4 u_3,
	\\
	\noalign{\medskip}
	[u_3,u_4]=  \gamma_5 u_1+\gamma_6 u_2 + (\gamma_1-\gamma_3) u_3 ,
\end{array}
\]  
where $\gamma_1$, $\dots$, $\gamma_6\in\mathbb{R}$. Since $\lambda_3\neq 0$, we consider the orthogonal basis $\hat u_i=\frac{1}{\lambda_3}u_i$ so that we can assume $\lambda_3=1$ working in the homothetic class of the initial metric. A direct  analysis of the components of the Bach tensor let us to clear some of the structure constants and completely determine the Bach-flat metrics. We proceed as follows. 

First, we compute 
$
\bach_{33}=\frac{1}{12}(
\gamma_6^2-8\gamma_5^2+16\gamma_2) 
$, which implies $\gamma_2= \frac{1}{16}( 8\gamma_5^2 - \gamma_6^2)$, and  as a consequence we get  
\[\small
\begin{array}{l}
	\phantom{..}
	-4\, \bach_{12} =  \gamma_6^2( 3\gamma_5 \gamma_6+2\gamma_1-3\gamma_3) ,
	
	\\
	\noalign{\medskip}
	
	-128\,\bach_{11} =   
	21 \gamma_6^4
	+144 \gamma_5^2 \gamma_6^2
	+ \gamma_5 \gamma_6 (64 \gamma_1  - 288 \gamma_3)
	-80 \gamma_1^2 + 240 \gamma_3^2
	- 160  \gamma_1 \gamma_3 + 96 \gamma_4 \gamma_6  ,
	
	\\
	\noalign{\medskip}
	
	\phantom{-}
	128\,\bach_{34} =
	63 \gamma_6^4
	+144 \gamma_5^2 \gamma_6^2
	+ \gamma_5 \gamma_6 (192 \gamma_1  - 192 \gamma_3 )
	+80 \gamma_1^2 + 80 \gamma_3^2
	-160 \gamma_1 \gamma_3 + 32 \gamma_4 \gamma_6  .
	
\end{array} 
\]
Note that if $\gamma_6\neq 0$ then $\gamma_3=\frac{1}{3} (3\gamma_5\gamma_6+2\gamma_1)$, and a direct calculation shows that
$\bach_{11} + 3\,\bach_{34} = \frac{1}{48}(40\gamma_1^2+63 \gamma_6^4)\neq 0$. 
Hence, necessarily $\gamma_6=0$  and, moreover, 
$\bach_{34}=\frac{5}{8}(\gamma_1-\gamma_3)^2$, which implies $\gamma_3=\gamma_1$.

At this point  we have  
$\gamma_2= \frac{1}{2} \gamma_5^2$,  $\gamma_6=0$ and    $\gamma_3=\gamma_1$, and a final straightforward calculation shows that the Bach tensor is determined by
\[
	\bach_{24} = \tfrac{3}{4}\gamma_1 (\gamma_1 \gamma_5+\gamma_4)
	\quad\text{and}\quad
	\bach_{44} = -\tfrac{3}{4}(\gamma_1 \gamma_5+\gamma_4) (3\gamma_1\gamma_5+\gamma_4).
\]
We conclude that $\gamma_4=-\gamma_1\gamma_5$ and thus the left-invariant metric corresponds to 
\[ 
	[u_1,u_4]= \gamma_1 u_1, \,\,\,
	[u_2,u_3]=  u_1, \,\,\,
	[u_2,u_4]= \tfrac{1}{2} \gamma_5^2 u_1 +\gamma_1 u_2 -\gamma_1\gamma_5 u_3, \,\,\,
	[u_3,u_4]=  \gamma_5 u_1.
\]  
Let $\zeta=-1$, $0$, $1$, depending on whether $\gamma_1$ is negative, zero,  or positive. Taking the basis
\[
	\tu_1=u_1,\quad
	\tu_2=u_2-\gamma_5 u_3,\quad
	\tu_3=u_3,\quad
	\tu_4=\gamma_5 u_2-\tfrac{1}{2}\gamma_5^2 u_3+u_4,
\]
if $\gamma_1=0$, and 
\[
	\tu_1=\tfrac{1}{|\gamma_1|^\frac{1}{2}} u_1,\,\,\,\,
	\tu_2=\tfrac{\zeta}{|\gamma_1|^\frac{1}{2}}\left(u_2-\gamma_5 u_3\right),\,\,\,\,
	\tu_3=\zeta u_3,\,\,\,\,
	\tu_4=\tfrac{1}{\gamma_1}\left(\gamma_5 u_2-\tfrac{1}{2}\gamma_5^2 u_3+u_4\right)
\]
if $\gamma_1\neq 0$, a direct calculations shows that 
 the Lie bracket transforms into
\[
	[\tu_1,\tu_4]= |\zeta|  \tu_1,\qquad
	[\tu_2,\tu_3]=  \tu_1, \qquad
	[\tu_2,\tu_4]=  |\zeta| \tu_2,
\]
while the inner product remains invariant if $\zeta=0$ and it is given by $\frac{1}{|\gamma_1|}\langle\cdot,\cdot\rangle$ if $\zeta\neq 0$. Since we are working at the homothetic level we can maintain the initial inner product in any case, remaining in the same homothetic class.
Now, a straightforward calculation shows that the  above metric  is not locally symmetric, locally conformally flat or Einstein. 
The Ricci operator of the above metric has eigenvalues $\{\zeta,\zeta,-\zeta,-\zeta\}$ and the geometry of the corresponding spacetime is completely different depending on $\zeta$ vanishes or not. 
If $\zeta\in\{-1,1\}$, the Ricci operator is diagonalizable  and the metric corresponds to that in Theorem~\ref{th:ce}-(D.i). Otherwise, if  $\zeta=0$, the Ricci operator is $2$-step nilpotent and a straightforward calculation shows that $\tu_4$ is a null parallel vector field   and   the curvature tensor satisfies $R(x,y)=0$ and $\nabla_xR=0$ for all $x,y\in \tu_4^\perp=\spn\{\tu_1,\tu_2,\tu_4\}$. Thus, the underlying structure is a left-invariant plane wave on the product $H_3\times\mathbb{R}$, and thus isomorphically homothetic to the Bach-flat metric discussed in Remark~\ref{re:producto}.

\smallskip

\subsubsection{Case $\lambda_1=0$,  $\lambda_2\neq 0$}\label{se:deg-2-2}

The Lie algebra structure is then given by
\[
\begin{array}{lll}
	[u_1,u_3]= \lambda_2 u_1, &
	[u_1,u_4]= \gamma_1 \lambda_2 u_1,  &
	[u_2,u_3]= \lambda_3 u_1, \quad
	\\
	\noalign{\medskip}
	[u_2,u_4]= (\gamma_1-\gamma_2)\lambda_3 u_1+ \gamma_2\lambda_2 u_2 , &
	[u_3,u_4]=  \gamma_3 u_1 + \gamma_4 u_2,
\end{array}
\]  
where $\gamma_1$, $\dots$, $\gamma_4\in\mathbb{R}$. Since $\lambda_2\neq 0$, we consider the orthogonal basis $\hat u_i=\frac{1}{\lambda_2}u_i$ so that we can assume $\lambda_2=1$ working in the homothetic class of the initial metric.  In this case, we compute a  Gr\"obner basis  of the ideal generated by the components of the Bach tensor, $\langle \bach_{ij}\rangle\subset \mathbb{R}[\lambda_3,\gamma_1$, $\gamma_2$, $\gamma_3$, $\gamma_4]$,  with respect to the   lexicographical order  (which consists of $48$ polynomials) and get  that the polynomials
\[
	\mathbf{g}_{1}= \gamma_4^5 (\gamma_3^2+\gamma_4^2)
	\quad\text{and}\quad
	\mathbf{g}_{2}= (20\gamma_3^2-7\gamma_4^2)(\gamma_3^2+\gamma_4^2)^2
\]
belong to the ideal. Thus, necessarily $\gamma_3=\gamma_4=0$.  Now, a direct calculation shows that $\bach_{33}=\frac{2}{3}(\gamma_1-\gamma_2)(2\lambda_3^4+3\lambda_3^2+1)$, which implies  $\gamma_2=\gamma_1$.  Under these assumptions, one easily checks that  the Bach tensor vanishes and the corresponding left-invariant metric is given by
\[ 
	[u_1,u_3]=   u_1, \quad
	[u_1,u_4]= \gamma_1   u_1,  \quad
	[u_2,u_3]= \lambda_3 u_1, \quad 
	[u_2,u_4]=   \gamma_1  u_2.
\]
Let $\{\tu_i\}$ be the basis given by the initial one ($\tu_i=u_i$)  if $\gamma_1=0$, or  
\[
	 \tu_1 	=  	  \tfrac{1}{|\gamma_1|^{\frac{1}{2}}} u_1 ,\quad
	 \tu_2	=  	  \tfrac{1}{|\gamma_1|^{\frac{1}{2}}} u_2 ,\quad
	 \tu_3  = 	u_3,\quad
	 \tu_4= \tfrac{1}{|\gamma_1|}  u_4 
\]
if $\gamma_1\neq 0$.
Let $\zeta=-1$, $0$, $1$, depending on whether $\gamma_1$ is negative, zero,  or positive. 
One easily checks that  the Lie bracket transforms into
\begin{equation}\label{eq:deg-2}
	[\tu_1,\tu_3]= \tu_1,\qquad
	[\tu_1,\tu_4]=  \zeta \tu_1, \qquad
	[\tu_2,\tu_3]=  \lambda_3  \tu_1,\qquad
	[\tu_2,\tu_4]= \zeta \tu_2,
\end{equation}
while  the inner product   obviously remains invariant if $\gamma_1 =0$ or,  otherwise, it  is given by $\frac{1}{|\gamma_1|} \langle\cdot,\cdot\rangle$. Since we are working at the homothetic level, we can maintain the initial inner product $\langle\cdot,\cdot\rangle$ in any case,  remaining in the same homothetic class. Now, a straightforward calculation shows that the above metric is not   locally conformally flat nor Einstein, and it  is locally symmetric if and only if $\lambda_3=\zeta=0$. Moreover, $(\tu_1,\tu_2,\tu_3,\tu_4)\mapsto (\tu_1,-\tu_2,\tu_3,\tu_4)$  defines an isometry interchanging $(\zeta,\lambda_3)$  and $(\zeta,-\lambda_3)$ and hence  we can restrict the parameter $\lambda_3$ to $\lambda_3\geq 0$. 
A direct calculation shows that the eigenvalues of the  Ricci operator of the above metric are given by 
$\{ \zeta(-1\pm\sqrt{\lambda_3^2+2}), \zeta(-2\pm\sqrt{\lambda_3^2+1}) \}$. 
If $\zeta\in\{-1,1\}$ then the four eigenvalues are different  and  the metrics correspond to Theorem~\ref{th:ce}-(D.ii).    If $\zeta=0$ then  the metric is never locally conformally flat nor locally symmetric. Moreover, the Ricci operator is $2$-step nilpotent and   $\tu_4$ is a null parallel vector field   so that     $R(x,y)=0$ and $\nabla_xR=0$ for all $x,y\in \tu_4^\perp=\spn\{\tu_1,\tu_2,\tu_4\}$. Thus, the underlying structure is a plane wave as in Theorem~\ref{th:plane-waves}-(D.ii).

\smallskip

\subsubsection{Case $\lambda_1\neq 0$}\label{se:deg-2-3}

In this case,  the Lie algebra structure becomes
\[
\begin{array}{l}
	[u_1,u_2]= \lambda_1 u_1, \quad
	[u_1,u_3]= \lambda_2 u_1, \quad
	[u_1,u_4]= \gamma_1 \lambda_1 u_1,  \quad
	[u_2,u_3]= \lambda_3 u_1, \quad
	\\
	\noalign{\medskip}
	[u_2,u_4]= \lambda_1\gamma_2 u_1 - \gamma_3\lambda_1\lambda_2 u_2 + \gamma_3 \lambda_1^2 u_3 , \quad
	\\
	\noalign{\medskip}
	[u_3,u_4]= -(\gamma_3\lambda_2\lambda_3-\gamma_2\lambda_2+(\gamma_1-\gamma_4)\lambda_3 ) u_1 -\gamma_4 \lambda_2 u_2 + \gamma_4 \lambda_1 u_3 , 
\end{array}
\]  
where $\gamma_1$, $\dots$, $\gamma_4\in\mathbb{R}$. As in the previous case, we use Gr\"obner bases. First, we compute a  Gr\"obner basis   of the ideal  $\langle \bar\bach_{ij}\rangle\subset \mathbb{R}[\lambda_1$, $\lambda_2$, $\lambda_3$, $\gamma_1$, $\gamma_2$, $\gamma_3$, $\gamma_4]$ with respect to the  graded reverse lexicographical order, where $\bar\bach_{ij}$ are obtained from the components $\bach_{ij}$ after simplifying the parameter $\lambda_1\neq0$ whenever possible. The Gröbner basis, consisting of $2183$ polynomials, is obtained after a long calculation (more than 12 hours with a last generation home computer) and contains the polynomial
\[
\mathbf{g}= \gamma_4^7 \lambda_1^5 (\lambda_2^2+\lambda_3^2) .
\]
Since $\lambda_1\neq 0$, it follows that either $\lambda_2=\lambda_3=0$ or $\gamma_4=0$.
If $\lambda_2=\lambda_3=0$ a direct calculation shows that $\bach_{11}=\frac{1}{6}\lambda_1^4$ and therefore  the metric cannot be Bach-flat. Thus, necessarily $\gamma_4=0$. Secondly, assuming  this last condition, we compute another Gr\"obner basis for   $\langle \bar\bach_{ij} \rangle\subset \mathbb{R}[\gamma_1$, $\gamma_2$, $\gamma_3$,   $\lambda_1$, $\lambda_2$, $\lambda_3]$ with respect to the  graded reverse lexicographical order  and get    the polynomials 
\[
\begin{array}{l}
	\mathbf{g}'_{1}= \lambda_1^4 (\lambda_1+2\gamma_1\lambda_2) (\lambda_1+2\lambda_2(\gamma_3\lambda_2+\gamma_1)),
	\\\noalign{\medskip}
	\mathbf{g}'_{2}= \lambda_1^2 (\lambda_1+2\gamma_1\lambda_2)^2 (2\gamma_2 \lambda_2^2-\lambda_1\lambda_3),
	\\\noalign{\medskip}
	\mathbf{g}'_{3}= \gamma_3\lambda_1^4 ((\lambda_1+2\gamma_1\lambda_2)(\lambda_1+2\gamma_2\lambda_3)
	+2\gamma_3\lambda_1(\lambda_2^2+\lambda_3^2)) ,
\end{array}
\]
among the $216$ polynomials of the basis.
The  expressions of $\mathbf{g}'_{1}$ and $\mathbf{g}'_{2}$ imply  that either $\lambda_1+2\gamma_1\lambda_2=0$ or, otherwise, $\lambda_1+2\lambda_2(\gamma_3\lambda_2+\gamma_1)=0$, $2\gamma_2 \lambda_2^2-\lambda_1\lambda_3=0$. Note that, in any case, $\lambda_1\neq 0$ implies $\lambda_2\neq 0$.  Next we show that the first condition implies the second one in the Bach-flat setting. Indeed, 
if $\lambda_1+2\gamma_1\lambda_2=0$, or equivalently $\gamma_1=-\frac{\lambda_1}{2\lambda_2}$, then 
$\mathbf{g}'_{3}= 2\gamma_3^2\lambda_1^6(\lambda_2^2+\lambda_3^2)
$. Hence   necessarily $\gamma_3=0$ and, in that case, a direct calculation shows that $\bach_{11}=-\frac{5(2\gamma_2 \lambda_2^2-\lambda_1\lambda_3)^4}{96\lambda_2^4}$. Thus  $2\gamma_2 \lambda_2^2-\lambda_1\lambda_3=0$ and, since $\gamma_3=0$, also $\lambda_1+2\lambda_2(\gamma_3\lambda_2+\gamma_1)=0$. As a consequence, we conclude that if $\gamma_4=0$ and the metric is Bach-flat, then $\lambda_1+2\lambda_2(\gamma_3\lambda_2+\gamma_1)=0$ and $2\gamma_2 \lambda_2^2-\lambda_1\lambda_3=0$, or equivalently $\gamma_1=-\frac{\lambda_1+2\gamma_3\lambda_2^2}{2\lambda_2}$ and  $\gamma_2=\frac{\lambda_1\lambda_3}{2\lambda_2^2}$. Finally, a straightforward calculation shows that, under these assumptions,  the Bach tensor vanishes and the corresponding left-invariant metric is given by
\[
\begin{array}{lll}
	[u_1,u_2]= \lambda_1 u_1, &
	[u_1,u_3]= \lambda_2 u_1, &
	[u_1,u_4]= -\frac{\lambda_1(\lambda_1+2\gamma_3\lambda_2^2)}{2\lambda_2}   u_1,  
	\\
	\noalign{\medskip}
	[u_2,u_3]= \lambda_3 u_1, &
	[u_3,u_4]= \frac{\lambda_1\lambda_3}{\lambda_2} u_1 &
	[u_2,u_4]=    \frac{\lambda_1^2\lambda_3}{2\lambda_2^2}  u_1 - \gamma_3\lambda_1\lambda_2 u_2 + \gamma_3 \lambda_1^2 u_3 .
\end{array}
\]
At this point, we introduce variables $\lambda_3' = \frac{\lambda_3}{\lambda_2}$ and $\gamma_3'=-\gamma_3\lambda_1 \lambda_2^2$. Moreover,  let $\zeta=-1$, $0$, $1$, depending on whether $\gamma_3'$ is negative, zero,  or positive. Considering the basis  $\{\tu_i\}$ defined by 
\[ 
	\tu_1
	= u_1 , \quad  
	\tu_2 =   u_2-\tfrac{\lambda_1}{\lambda_2} u_3  , \quad
	\tu_3=
	\tfrac{1}{\lambda_2} u_3, \quad
	\tu_4=   \lambda_1 u_2-\tfrac{\lambda_1^2}{2\lambda_2}u_3 +\lambda_2 u_4 	 
\]
if $\gamma_3'=0$, and
\[ 
\begin{array}{ll}
	\tu_1
	= \tfrac{1}{|\gamma_3'|^{\frac{1}{2}}}  u_1 ,    &
	\quad
	\tu_2 = \tfrac{1}{|\gamma_3'|^{\frac{1}{2}}}   \left(u_2-\tfrac{\lambda_1}{\lambda_2} u_3\right) , 
	\\
	\noalign{\medskip}
	\tu_3= 	\tfrac{1}{\lambda_2} u_3,   &
	\quad 
	\tu_4=  \tfrac{1}{|\gamma_3'|}   \left(\lambda_1 u_2-\tfrac{\lambda_1^2}{2\lambda_2}u_3 +\lambda_2 u_4\right)	   
\end{array}
\]
if $\gamma_3'\neq 0$, 
a direct calculation shows that  the Lie bracket transforms into
\[
[\tu_1,\tu_3]= \tu_1,\qquad
[\tu_1,\tu_4]=  \zeta \tu_1, \qquad
[\tu_2,\tu_3]=  \lambda_3'  \tu_1,\qquad
[\tu_2,\tu_4]= \zeta \tu_2,
\]
while  the inner product   remains invariant if $\zeta =0$ or, otherwise, it  is given by $\frac{1}{|\gamma_3'|} \langle\cdot,\cdot\rangle$. Since we are working at the homothetic level, we can maintain the initial inner product $\langle\cdot,\cdot\rangle$ in any case,  remaining in the same homothetic class. Hence,  the setting  is exactly the same as in $\S$\ref{se:deg-2-2} (see Equation~\eqref{eq:deg-2}) and therefore the above metric is a plane wave if $\zeta=0$ (corresponding to Theorem~\ref{th:plane-waves}-(D.ii)) while, otherwise,  it corresponds to   Theorem~\ref{th:ce}--(D.ii).

\section{Semi-direct extensions with Riemannian normal subgroup $H_3$} \label{se:Riemann}

In this section we consider   left-invariant  Lorentzian metrics which are obtained as extensions 
of the three-dimensional unimodular Riemannian Heisenberg Lie group $H_3$. 
Left-invariant Riemannian metrics on three-dimensional unimodular Lie groups were described by Milnor \cite{Milnor} using the self-dual structure tensor $L$ given by $L(X\times Y)=[X,Y]$, where ``$\times$" denotes the vector-cross product $\langle X\times Y,Z\rangle=\operatorname{det}(X,Y,Z)$. Self-duality of $L$ ensures the existence of an orthonormal basis of $\mathfrak{h}_3$ diagonalizing the structure tensor in the positive definite case. Therefore any Riemannian left-invariant metric on $H_3$ is homothetic to the one determined at the Lie algebra level by 
an orthonormal basis $\{\mathbf{v}_1, \mathbf{v}_2, \mathbf{v}_3\}$ of $\mathfrak{h}_3$ such that (see \cite{Milnor})
$$
[\mathbf{v}_3,\mathbf{v}_2]=0,
\quad 
[\mathbf{v}_3,\mathbf{v}_1]=0, 
\quad 
[\mathbf{v}_1,\mathbf{v}_2]= \lambda \mathbf{v}_3,\qquad \lambda\neq 0.
$$
The algebra of all derivations of $\mathfrak{h}_3$ is given, with respect to the basis $\{ \mathbf{v}_1, \mathbf{v}_2, \mathbf{v}_3\}$, by
\begin{equation}\label{eq:EGR-derivation}
\operatorname{der}(\mathfrak{h}_3)=\left \{  \left( \begin{array}{ccc}
\alpha_{11} & \alpha_{12} & 0 
\\
\alpha_{21} & \alpha_{22} & 0 
\\
\alpha_{31} & \alpha_{32} &\alpha_{11}+\alpha_{22}
\end{array} \right) ;\, \alpha_{ij} \in \mathbb{R} \right \}.
\end{equation} 
For any semi-direct extension $\mathfrak{h}_3\rtimes\mathfrak{r}$, there is a basis $\{\mathbf{v}_1, \mathbf{v}_2, \mathbf{v}_3,\mathbf{v}_4\}$ so that $\mathfrak{g}=\mathfrak{h}_3\rtimes\mathfrak{r}=\operatorname{span}\{\mathbf{v}_1, \mathbf{v}_2, \mathbf{v}_3\}\oplus\mathbb{R}\mathbf{v}_4$.
Since $\mathbb{R}\mathbf{v}_4$ is timelike but not necessarily orthogonal to $\mathfrak{h}_3$, we set $\bar e_4= \mathbf{v}_4-\sum_i \langle \mathbf{v}_4,\mathbf{v}_i \rangle\mathbf{v}_i$ and normalize it to get an orthonormal basis $\{ e_1,\dots,e_4\}$ of $\mathfrak{g}= \mathfrak{h}_3\oplus\mathbb{R}$ with $e_4$ timelike, where $e_i=\mathbf{v}_i$ $(i=1,2,3)$, so that 
\[
\begin{array}{ll}
	[e_1,e_2]=\lambda e_3,  &
	[e_1,e_4]=\gamma_1 e_1+\gamma_2 e_2 + \gamma_3 e_3,  
	\\
	\noalign{\medskip}
	[e_2,e_4]=\gamma_4 e_1 + \gamma_5  e_2 + \gamma_6 e_3,  &
	[e_3,e_4]=(\gamma_1 +\gamma_5)e_3,
\end{array}
\]
where $\lambda\neq 0$ and $\gamma_1$, $\dots$, $\gamma_6\in\mathbb{R}$. 

In what follows  we determine the Bach-flat metrics analyzing  two cases by separate, depending on whether $\gamma_5$ vanishes or not.

\subsection{\bf Case $\boldsymbol{\gamma_5=0}$}\label{se:Riemann-1}

Since $\lambda\neq 0$, we can assume $\lambda=1$ working in the homothetic class of the initial metric, just taking the orthogonal basis $\hat e_i=\frac{1}{\lambda}e_i$. 
Considering   the ideal  $\langle \bach_{ij}\rangle\subset \mathbb{R}[ \gamma_2$, $\gamma_3$, $\gamma_4$,  $\gamma_6$, $\gamma_1]$,  where $\bach_{ij}=\bach(e_i,e_j)$ are the polynomials given by the components of the Bach tensor in the orthonormal basis $\{ e_i\}$, and fixing the   lexicographical order we make use of     Gr\"obner bases. Thus we get 50 polynomials, being one of them the polynomial 
\[
\mathbf{g}= \gamma_1^4 \gamma_6^2 (\gamma_1^2+1) ( \gamma_1^2+2) ( \gamma_1^2+4) 
(8 \gamma_1^2+5) ( 16 \gamma_1^2+25) (25 \gamma_1^2+1) .
\]
Hence, either $\gamma_1=0$ or $\gamma_1\neq 0$ and $\gamma_6=0$. Moreover, if $\gamma_1=0$, we get that   
\[
\mathbf{g}'_1= \gamma_4^4
\quad\text{and}\quad
\mathbf{g}'_2= \gamma_2^2-8\gamma_3^2+\gamma_4^2-8\gamma_6^2+8
\]
are polynomials in the ideal $\langle \bach_{ij}\rangle$, and the same occur with the polynomials
\[
\mathbf{g}''_1= \gamma_3\gamma_4^4
\quad\text{and}\quad
\mathbf{g}''_2= \gamma_1\gamma_2\gamma_3 (\gamma_1^2+1)
\]
if   $\gamma_6=0$. Next we examine the different possibilities separately.

\smallskip

\subsubsection{Case $\gamma_1=\gamma_4=0$, $\gamma_2^2-8\gamma_3^2 -8\gamma_6^2+8=0$}\label{se:Riemann-1-1}

In this case, a direct calculation shows that $\bach_{11}=\frac{9}{128}\gamma_2^2(7\gamma_2^2+16\gamma_3^2)$, which implies $\gamma_2=0$, and  the Bach tensor vanishes under this last assumption.  Thus, the corresponding  left-invariant metric is given by
\[ 
	[e_1,e_2]=  e_3,  \qquad
	[e_1,e_4]=   \gamma_3 e_3, \qquad  
	[e_2,e_4]=   \varepsilon \sqrt{1-\gamma_3^2}\, e_3,
\]
where $\gamma_3\in[-1,1]$ and $\varepsilon^2=1$.
The Ricci operator of the above metric is $2$-step nilpotent and   
$\xi=  \varepsilon \sqrt{1-\gamma_3^2}\,   e_1 -\gamma_3  e_2 +   e_4$ 
is a null parallel vector field. Moreover, the curvature tensor satisfies     $R(x,y)=0$ and $\nabla_xR=0$ for all $x,y\in \xi^\perp=\spn\{\xi,e_3,\theta\}$, where 
\[
\theta=\left\{
\begin{array}{ll}
	e_2 & \text{if $\gamma_3=0$},
	\\
	\noalign{\medskip}
	e_1 + \varepsilon\sqrt{1-\gamma_3^2}\, e_4 & \text{if $\gamma_3\neq 0$}.
\end{array}
\right.
\]
Thus, the underlying structure is a plane wave. 
Moreover, considering the basis
	$$
	\te_1=e_1,\quad \te_2=e_2, \quad \te_3=e_3, \quad 
 \te_4=\varepsilon\sqrt{1-\gamma_3^2}  e_1-\gamma_3e_2+e_4 ,
	$$
	one has that the only non-zero bracket is $[\te_1,\te_2]=\te_3$, which shows that the underlying Lie group is the product $H_3\times\mathbb{R}$.	
	It is therefore isomorphically homothetic to the Bach-flat metric discussed in Remark~\ref{re:producto}.

\smallskip

\subsubsection{Case $\gamma_1\neq 0$, $\gamma_3=\gamma_6=0$}\label{se:Riemann-1-2}

By a straightforward  calculation one gets the following components of the Bach tensor
\[\small
\begin{array}{l}
	-24\, \bach_{33}= 
	(\gamma_2 + \gamma_4)^2 ( 
	4 (\gamma_2 + \gamma_4)^2 - 12  \gamma_2 \gamma_4 + 
	3) - \gamma_1^2 ((\gamma_2 + \gamma_4)^2 + 
	12 \gamma_2 \gamma_4 - 28  )  + 20,
	\\
	\noalign{\medskip}
	-24\, \bach_{44}= 
	(\gamma_2 + \gamma_4)^2 (  
	12 (\gamma_2 + \gamma_4)^2 - 36 \gamma_2 \gamma_4 + 
	1 ) + \gamma_1^2 (  
	13 (\gamma_2 + \gamma_4)^2 - 36 \gamma_2 \gamma_4 - 4) - 4,
\end{array}
\]
which implies 
\[  
\begin{array}{l}
-12 (\bach_{33} +3\,\bach_{44}) =  
(\gamma_2 + \gamma_4)^2 (20 (\gamma_2 + \gamma_4)^2 - 
60 \gamma_2 \gamma_4 + 
3) 
\\
\noalign{\medskip}
\phantom{-12 (\bach_{33} +3\,\bach_{44}) =}
+ \gamma_1^2 (19 (\gamma_2 + \gamma_4)^2 - 
60  \gamma_2 \gamma_4 + 8 ) + 4.
\end{array}
\]
Clearly this last expression is always strictly positive and therefore  the metric cannot be Bach-flat in this case.

\smallskip

\subsubsection{Case $\gamma_1\gamma_3\neq 0$, $\gamma_2=\gamma_4=\gamma_6=0$}\label{se:Riemann-1-3}

A direct calculation shows that, in this case, $\bach_{22}=-\frac{1}{6}(3\gamma_1^2+\gamma_3^2+3)(\gamma_3^2-1)$. Hence $\gamma_3=\varepsilon$, with $\varepsilon^2=1$, and a straightforward calculation shows that the Bach tensor vanishes. Thus, the left-invariant metric corresponds to
\begin{equation}\label{eq:riemann-1}
	[e_1,e_2]=  e_3,  \qquad
	[e_1,e_4]=\gamma_1 e_1 + \varepsilon e_3,   \qquad 
	[e_3,e_4]= \gamma_1  e_3, 
\end{equation}
which is not locally symmetric, locally conformally flat nor Einstein. 
Moreover, $(e_1,e_2,e_3,e_4)\mapsto (-e_1,-e_2,e_3,-e_4)$ and $(e_1,e_2,e_3,e_4)\mapsto (e_1,e_2,e_3,-e_4)$  define   isometries interchanging $(\varepsilon,\gamma_1)$ with $(\varepsilon,-\gamma_1)$  and $(\varepsilon,\gamma_1)$ with $(-\varepsilon,-\gamma_1)$, respectively,  and hence  we can restrict the parameters to $\varepsilon=1$ and  $\gamma_1> 0$. This corresponds to the metric in Theorem~\ref{th:ce}--(R).

\smallskip

\subsection{\bf Case $\boldsymbol{\gamma_5\neq 0}$}\label{se:Riemann-2}

Taking   the orthogonal basis $\hat e_i=\frac{1}{\gamma_5}e_i$ we  can assume $\gamma_5=1$ continuing  in the homothetic class of the initial metric.  We fix the   polynomial ring $\mathbb{R}[ \gamma_2$, $\gamma_3$, $\gamma_4$, $\gamma_6$, $\gamma_1$, $\lambda]$. 
Considering the ideal $\langle \overline{\bach}_{ij}\rangle$, where $\overline{\bach}_{ij}$ are the components of the Bach tensor where we have simplified the parameter $\lambda\neq 0$ whenever possible, and computing a   Gr\"obner basis with respect to the graded reverse lexicographical order (consisting of 211 polynomials) we  get that the polynomial 
\[
\mathbf{g}_1=
((12 \gamma_2 + 5  \gamma_4) \gamma_3^2 - (5 \gamma_2 + 
12 \gamma_4 ) \gamma_6^2 - 
17 ( \gamma_1 - 1) \gamma_3 \gamma_6) 
(\lambda^2 + (\gamma_1 + 1)^2) \lambda^{10}
\]
belongs to the ideal $\langle \bar\bach_{ij}\rangle$. Hence,  since $\lambda\neq 0$, it follows that the polynomial     
\begin{equation}\label{eq:Riemann-2-1}
\mathbf{\tilde g}_1=
	(12 \gamma_2 + 5  \gamma_4) \gamma_3^2 
	- (5 \gamma_2 + 12 \gamma_4 ) \gamma_6^2 
	-  17 ( \gamma_1 - 1) \gamma_3 \gamma_6
\end{equation}
must vanish. 
Now, for the ideal $\langle \overline{\bach}_{ij}\cup \{\mathbf{\tilde g}_1\}\rangle$, we consider the   lexicographical order and obtain, after computing a Gr\"obner basis, that
\[
\begin{array}{l}
	\mathbf{g}_2 = 
	-\gamma_6 
	(\lambda^2 - (\gamma_1 + 1) \gamma_6^2) 
	(\lambda^2 + (\gamma_1 + 1)^2) 
	(\lambda^2 + 25 (\gamma_1 + 1)^2) 
	(2 \lambda^2 + (\gamma_1 + 1)^2) 
	\\
	\noalign{\medskip}
	\phantom{\mathbf{g_2} = }
	\times
	(4 \lambda^2 + (\gamma_1 + 1)^2) 
	(5 \lambda^2 + 8 (\gamma_1 + 1)^2) 
	(25 \lambda^2 + 16 (\gamma_1 + 1)^2)
\end{array}
\]
belongs to the ideal. Thus, either $\gamma_6=0$ or $\gamma_6\neq0$ and $\lambda^2 - (\gamma_1 + 1) \gamma_6^2=0$.  Next we examine these two  cases separately.

\smallskip

\subsubsection{Case $\gamma_6=0$}\label{se:Riemann-2-1}

In this case  we show that  there is no Bach-flat metrics. We make use of Gr\"obner bases  considering  the lexicographical order in all the cases. First, we consider the polynomial ring $\mathbb{R}[ \gamma_2$, $\gamma_3$, $\gamma_4$,   $\gamma_1$, $\lambda]$. Taking  the ideal   $\langle \overline{\bach}_{ij}\rangle$   we   get,  after computing a Gröbner basis, that
\[
\mathbf{g}'_1=
\gamma_3 \lambda^6 
(\lambda^2 + 4) 
(5 \lambda^2 + 2) 
(16 \lambda^2 + 7225) 
(32 \lambda^2 + 9) 
(256 \lambda^2 + 169)
\]
belongs to $\langle \bar\bach_{ij}\rangle$, which implies $\gamma_3=0$. Adding this  new condition, we repeat the process in the polynomial ring $\mathbb{R}[ \gamma_2$,   $\gamma_4$,   $\gamma_1$, $\lambda]$ to obtain the polynomial 
\[
\begin{array}{l}
	\mathbf{g}'_2=
	-\lambda^2 
	(7 \lambda^2 - 2 (\gamma_1 + 1)^2) 
	(\lambda^2 + (\gamma_1 + 1)^2) 
	(\lambda^2 + 2 (\gamma_1 + 1)^2) 
	\\
	\noalign{\medskip}
	\phantom{\mathbf{g}'_2=}
	\times
	(4 \lambda^2 + (\gamma_1 + 1)^2) 
	(5 \lambda^2 + 8 (\gamma_1 + 1)^2).
\end{array}
\]
Thus,   $\mathbf{\tilde g}'_2=7 \lambda^2 - 2 (\gamma_1 + 1)^2$ must vanish. Finally, we consider the polynomials in the ring $\mathbb{R}[ \gamma_1$,  $\lambda$,   $\gamma_2$, $\gamma_4]$ and 
 compute a last  Gr\"obner basis for the ideal $\langle \overline{\bach}_{ij}\cup\{\mathbf{\tilde g}'_2\}\rangle$. As a result, we get that  
\[
\begin{array}{l}
	\mathbf{g}'_3=
	(\gamma_4^2 + 1)^4 
	(441 \gamma_2^4 
	+ 2 \gamma_2^2 (641 \gamma_4^2 + 770) 
	+ (7 \gamma_4^2 + 22) (63 \gamma_4^2 + 22) 
	\\
	\noalign{\medskip}
	\phantom{\mathbf{g}'_3=(\gamma_4^2 + 1)^4}
	- 24 \gamma_2 \gamma_4 (35 \gamma_2^2 + 35 \gamma_4^2 + 88))
\end{array}
\]
belongs to the ideal. Since $\mathbf{g}'_3$ does not vanish, we conclude that there is no Bach-flat metrics in this case.

\smallskip

\subsubsection{Case $\gamma_6\neq0$, $\lambda^2 - (\gamma_1 + 1) \gamma_6^2=0$}\label{se:Riemann-2-2}

To solve this case we work in the polynomial ring $\mathbb{R}[ \gamma_1$, $\gamma_6$, $\gamma_4$,   $\gamma_2$, $\gamma_3$, $\lambda]$ and use the lexicographical order for computing Gr\"obner bases. We start with the ideal $\langle\overline{\bach}_{ij}\cup\{\lambda^2 - (\gamma_1 + 1) \gamma_6^2\}\rangle$ and get that 
\[
\begin{array}{l}
	\mathbf{g}''_1 = 
		-\lambda^6
		(\lambda^2 + 1) 
		(4 \lambda^2 + 9) 
		(25 \lambda^2 + 4) 
		(45 \lambda^2 + 98) 
		(225 \lambda^2 + 256) 
		(256 \lambda^2 + 1369) 
		\\
		\noalign{\medskip}
		\phantom{\mathbf{g}''_1 = }
		\times
		(968 \lambda^2 + 1521)  
		(2304 \lambda^4 + 46441425 \lambda^2 + 1028805625) 
		\\
		\noalign{\medskip}
		\phantom{\mathbf{g}''_1 = }
		\times
		(3211264 \lambda^4 + 34731953 \lambda^2 + 35796289)
		(\gamma_2^2 \lambda^2 - (\gamma_2^2 + 1) \gamma_3^2)
\end{array}
\]
belongs to the ideal. Since $\lambda\neq 0$, only the last factor, $\mathbf{\tilde g}''_1  = \gamma_2^2 \lambda^2 - (\gamma_2^2 + 1) \gamma_3^2$, may vanish. Finally, we use again the polynomial $\mathbf{\tilde g}_1$ given by Equation~\eqref{eq:Riemann-2-1}  and considering the ideal $\langle\overline{\bach}_{ij}\cup\{\lambda^2 - (\gamma_1 + 1) \gamma_6^2$, $\mathbf{\tilde g}''_1$, $\mathbf{\tilde g}_1\}\rangle$ we obtain that the polynomials
\[ 
	\mathbf{g}''_2 = (\gamma_4 - \gamma_2) \lambda^4, \quad
	\mathbf{g}''_3 = ( \gamma_2 \gamma_6 - \gamma_3 ) \lambda^4
	\quad\text{and}\quad 
	\mathbf{g}''_4 = (79 \gamma_1 + 11 \gamma_2^2 - 	90 \gamma_2 \gamma_4) \lambda^2 
\]
belong to the ideal.
These expressions, together with $\lambda^2 - (\gamma_1 + 1) \gamma_6^2=0$, imply
\[
	\gamma_4=\gamma_2,\quad
	\gamma_3=\gamma_2 \gamma_6,\quad
	\gamma_1=\gamma_2^2,\quad
	\lambda=\varepsilon\gamma_6\sqrt{\gamma_2^2+1},
\]
where $\varepsilon^2=1$, and a straightforward calculation shows that, under these conditions, the Bach tensor vanishes. Hence, the left-invariant  metric is determined by 
\[
\begin{array}{ll}
	[e_1,e_2]=\varepsilon\gamma_6\sqrt{\gamma_2^2+1} e_3,  &
	[e_1,e_4]=\gamma_2^2 e_1+\gamma_2 e_2 + \gamma_2 \gamma_6 e_3,  
	\\
	\noalign{\medskip}
	[e_2,e_4]=\gamma_2 e_1 +   e_2 + \gamma_6 e_3,  &
	[e_3,e_4]=(\gamma_2^2 + 1)e_3.
\end{array}
\] 
Let $\gamma_6' = \frac{\sqrt{\gamma_2^2+1}}{\gamma_6}$ and consider a new basis  $\{\te_i\}$ defined by 
\[ 
\begin{array}{ll}
	\te_1
	= \frac{\gamma_6'}{(\gamma_2^2+1)^\frac{3}{2}} \left(  \gamma_2 e_1 + e_2\right) ,   &
	\quad
	\te_2 = -\frac{\gamma_6'}{(\gamma_2^2+1)^\frac{3}{2}} \left(  e_1 - \gamma_2 e_2\right) , 
	\\
	\noalign{\medskip}
	\te_3= 	\frac{\varepsilon\gamma_6'}{\gamma_2^2+1} e_3 ,  &
	\quad 
	\te_4=  \frac{\gamma_6'}{\gamma_2^2+1} e_4. 	   
\end{array}
\]
Now, a direct calculation shows that  the Lie bracket transforms into
\[
	[\te_1,\te_2]=  \te_3,  \qquad
	[\te_1,\te_4]=\gamma_6' \te_1 + \varepsilon \te_3,   \qquad 
	[\te_3,\te_4]= \gamma_6'  \te_3, 
\]
while  the inner product becomes $\left(\frac{\gamma_6'}{\gamma_2^2+1}\right)^2 \langle\cdot,\cdot\rangle$. Since we are working at the homothetic level, we can maintain the initial inner product   remaining in the same homothetic class. Thus,  we recover the case in $\S$\ref{se:Riemann-1-3} (see Equation~\eqref{eq:riemann-1}) and therefore the above metric  corresponds to the situation  given by   Theorem~\ref{th:ce}-(R).

\section{Semi-direct extensions with Lorentzian normal subgroup $H_3$}\label{se:EGR-Lorentzian}

In this section we deal with  left-invariant  Lorentzian metrics whose restriction to the three-dimensional unimodular Lie group $H_3$ is of Lorentzian signature.   
Left-invariant Lorentzian metrics on the Heisenberg group were described in \cite{Rahmani} by using Milnor type frames.
Let $L(X\times Y)=[X,Y]$, where ``$\times$" denotes the vector-cross product $\langle X\times Y,Z\rangle=\operatorname{det}(X,Y,Z)$, be the structure tensor. Self-duality of $L$ holds true in the unimodular case as well as in the Riemannian case, but due to the Lorentzian signature,  $L$ may have non-trivial Jordan normal form.
Since $L$ must have eigenvalues $\{ 0,0,\lambda_3\}$ one has that the only possible Jordan normal forms are as follows:
\begin{enumerate}
\item[\bf{Ia.}] $L$ is real diagonalizable. Hence there exists an orthonormal basis $\left\{e_1,e_2,e_3\right\}$, where we assume $e_3$ to be timelike, so that $L(e_i)=\lambda_i e_i$ with $\lambda_1=\lambda_2=0$.
\item[\bf{II.}] $L$ has a double root of its minimal polynomial. Then it is two-step nilpotent and there exists a pseudo-orthonormal basis  $\left\{u_1,u_2,u_3\right\}$ so that
\[
L=\begin{pmatrix}
0&0&0\\
\pm 1 &0&0\\
0&0&0
\end{pmatrix},
\quad
\text{where}\quad 
\langle u_1,u_2\rangle=\langle u_3,u_3\rangle=1 \,.
\]
\end{enumerate}
Therefore inner products on $\mathfrak{h}_3\rtimes\mathfrak{r}$ are described as in Section~\ref{se:Riemann} 	
	by considering their restriction to the subalgebra $\mathfrak{h}_3$, which is of type Ia (distinguishing the two possibilities whether the kernel of the structure operator is positive definite or Lorentzian) or two-step nilpotent as in II.
 Next we analyze the vanishing of the Bach tensor in those three cases separately.

\subsection{The structure operator is diagonalizable of rank one with positive definite kernel}\label{se:Lorentz-1}

In this case, using \eqref{eq:EGR-derivation} and proceeding as in Section~\ref{se:Riemann}, one gets that there exists an orthonormal basis $\{e_1,e_2,e_3,e_4\}$ of $\mathfrak{g}=\mathfrak{h}_3\rtimes \mathfrak{r}$, with $e_3$ timelike,  where  $\mathfrak{h}_3=\spn\{e_1,e_2,e_3\}$ and $\mathfrak{r}=\spn\{e_4\}$, so that 
\[
\begin{array}{ll}
	[e_1,e_2]=-\lambda e_3,  &
	[e_1,e_4]=\gamma_1 e_1+\gamma_2 e_2 + \gamma_3 e_3,  
	\\
	\noalign{\medskip}
	[e_2,e_4]=\gamma_4 e_1 + \gamma_5  e_2 + \gamma_6 e_3,  &
	[e_3,e_4]=(\gamma_1 +\gamma_5)e_3,
\end{array}
\] 
where $\lambda\neq 0$ and $\gamma_1$, $\dots$, $\gamma_6\in\mathbb{R}$. The analysis of this case is   analogous to the one carried out in Section~\ref{se:Riemann}. However,  the different behavior of the restriction of the inner product to the subalgebra $\mathfrak{h}_3$ (Lorentzian instead of Riemannian) is crucial and leads to the non-existence of Bach-flat metrics. Next we schematize the process distinguishing
the cases $\gamma_5=0$ and $\gamma_5\neq 0$, but omitting the details which coincide with the study developed in Section~\ref{se:Riemann}.

\smallskip

\subsubsection{Case $\gamma_5=0$}\label{se:Lorentz-1-1}

Since $\lambda\neq 0$, we can assume $\lambda=1$ and work in the homothetic class of the initial metric.
Proceeding as in Section~\ref{se:Riemann-1} we 
consider the ideal generated by the polynomials $\bach_{ij}\in\mathbb{R}[\gamma_2,\gamma_3,\gamma_4,\gamma_6,\gamma_1]$ and compute a Gröbner basis (consisting of 50 polynomials) with respect to the lexicographical order to
get that 
\[
\mathbf{g}= \gamma_1^4 \gamma_6^2 (\gamma_1^2+1) ( \gamma_1^2+2) ( \gamma_1^2+4) 
(8 \gamma_1^2+5) ( 16 \gamma_1^2+25) (25 \gamma_1^2+1) 
\]
belongs to the ideal $\langle\bach_{ij}\rangle$.  Hence, either $\gamma_1=0$ or $\gamma_1\neq 0$, $\gamma_6=0$. If $\gamma_1=0$, then   
\[
\mathbf{g}'_1=  \gamma_2^2+8\gamma_3^2+\gamma_4^2+8\gamma_6^2+8
\]
is a  polynomial in $\langle \bach_{ij}\rangle$ and therefore there is no Bach-flat metric in this case. Moreover, if $\gamma_6=0$,   we obtain the polynomials 
\[
\mathbf{g}''_1= \gamma_3\gamma_4^4
\quad\text{and}\quad
\mathbf{g}''_2= \gamma_1\gamma_2\gamma_3 (\gamma_1^2+1).
\]
Now, if $\gamma_3=0$,   proceeding as in Section~\ref{se:Riemann-1-2} we get
\[  
\begin{array}{l}
	12 (\bach_{33} +3\,\bach_{44}) =  
	(\gamma_2 + \gamma_4)^2 (20 (\gamma_2 + \gamma_4)^2 - 
	60 \gamma_2 \gamma_4 + 
	3) 
	\\
	\noalign{\medskip}
	\phantom{-12 (\bach_{33} +3\,\bach_{44}) =}
	+ \gamma_1^2 (19 (\gamma_2 + \gamma_4)^2 - 
	60  \gamma_2 \gamma_4 + 8 ) + 4,
\end{array}
\]
while if $\gamma_3\neq0$ then necessarily $\gamma_2=\gamma_4=0$ and a direct calculation shows that 
$\bach_{11}=\frac{1}{6} (\gamma_3^2+1) (5 \gamma_1^2 + 3 \gamma_3^2+3)$. Thus, in any case, the Bach tensor does not vanish.

\smallskip

\subsubsection{Case $\gamma_5\neq 0$}\label{se:Lorentz-1-2}

Taking   the orthogonal basis $\hat e_i=\frac{1}{\gamma_5}e_i$ we  can assume $\gamma_5=~1$.  
Consider the ideal generated by the polynomials $\bar\bach_{ij}\in\mathbb{R}[\gamma_2,\gamma_3,\gamma_4,\gamma_6,\gamma_1,\lambda]$, where $\bar\bach_{ij}$ are obtained from the components $\bach_{ij}$ after simplifying $\lambda\neq0$ whenever possible.
We proceed as in Section~\ref{se:Riemann-2} to construct a Gröbner basis (which consists of 211 polynomials)  with respect to the graded reverse lexicographical order and get that   
\[
\mathbf{g}_1=
((12 \gamma_2 + 5  \gamma_4) \gamma_3^2 - (5 \gamma_2 + 
12 \gamma_4 ) \gamma_6^2 - 
17 ( \gamma_1 - 1) \gamma_3 \gamma_6) 
(\lambda^2 + (\gamma_1 + 1)^2) \lambda^{10}
\]
is a polynomial in the ideal  $\langle \bar\bach_{ij}\rangle$.  Therefore  the polynomial     
\begin{equation}\label{eq:Lorentz-1-2-1}
	\mathbf{\tilde g}_1=
	(12 \gamma_2 + 5  \gamma_4) \gamma_3^2 
	- (5 \gamma_2 + 12 \gamma_4 ) \gamma_6^2 
	-  17 ( \gamma_1 - 1) \gamma_3 \gamma_6
\end{equation}
must vanish. 
Now, considering  the ideal $\langle  \bar\bach_{ij}\cup \{\mathbf{\tilde g}_1\}\rangle$ and the lexicographical order, we obtain that the polynomial
\[
\begin{array}{l}
	\mathbf{g}_2 = 
	\gamma_6 
	(\lambda^2 + (\gamma_1 + 1) \gamma_6^2) 
	(\lambda^2 + (\gamma_1 + 1)^2) 
	(\lambda^2 + 25 (\gamma_1 + 1)^2) 
	(2 \lambda^2 + (\gamma_1 + 1)^2) 
	\\
	\noalign{\medskip}
	\phantom{\mathbf{g_2} = }
	\times
	(4 \lambda^2 + (\gamma_1 + 1)^2) 
	(5 \lambda^2 + 8 (\gamma_1 + 1)^2) 
	(25 \lambda^2 + 16 (\gamma_1 + 1)^2)
\end{array}
\]
belongs to the ideal.
Thus, either $\gamma_6=0$ or $\gamma_6\neq0$, $\lambda^2 + (\gamma_1 + 1) \gamma_6^2=0$.  
If $\gamma_6=0$, proceeding   as in Section~\ref{se:Riemann-2-1} we get exactly the same polynomials $\mathbf{g}'_1$, $\mathbf{g}'_2$ and $\mathbf{g}'_3$, so we conclude that there is no Bach-flat metric in this case. Now, if $\gamma_6\neq0$ and  $\lambda^2 + (\gamma_1 + 1) \gamma_6^2=0$, we proceed as in Section~\ref{se:Riemann-2-2} to get the polynomial 
\[
\begin{array}{l}
	\mathbf{g}''_1 = 
	\lambda^6
	(\lambda^2 + 1) 
	(4 \lambda^2 + 9) 
	(25 \lambda^2 + 4) 
	(45 \lambda^2 + 98) 
	(225 \lambda^2 + 256) 
	(256 \lambda^2 + 1369) 
	\\
	\noalign{\medskip}
	\phantom{\mathbf{g}''_1 = }
	\times
	(968 \lambda^2 + 1521)  
	(2304 \lambda^4 + 46441425 \lambda^2 + 1028805625) 
	\\
	\noalign{\medskip}
	\phantom{\mathbf{g}''_1 = }
	\times
	(3211264 \lambda^4 + 34731953 \lambda^2 + 35796289)
	(\gamma_2^2 \lambda^2 + (\gamma_2^2 + 1) \gamma_3^2)
\end{array}
\]
in the ideal $\langle\bar\bach_{ij}\cup\{\lambda^2 + (\gamma_1 + 1) \gamma_6^2\}\rangle$. Hence, $\mathbf{\tilde g}''_1  = \gamma_2^2 \lambda^2 + (\gamma_2^2 + 1) \gamma_3^2$ must vanish. Thus,   $\gamma_2=\gamma_3=0$ and the polynomial $\mathbf{\tilde g}_1$ given in Equation~\eqref{eq:Lorentz-1-2-1} reduces to $\mathbf{\tilde g}_1 = -12\gamma_4\gamma_6^2$, which implies $\gamma_4=0$ since $\gamma_6\neq 0$.  
Finally, using $\gamma_2=\gamma_3=\gamma_4=0$ and $\lambda^2 + (\gamma_1 + 1) \gamma_6^2=0$, a straightforward calculation shows that 
\[ 
\begin{array}{l}
	\frac{24\gamma_6^6}{\lambda^2+\gamma_6^2}\, \bach_{11}= 
	4(5 \gamma_6^2 + 4) \lambda^4
	+3 (4 \gamma_6^4 - 3 \gamma_6^2 + 16 )  \gamma_6^2 \lambda^2
	-(4 \gamma_6^4 - 11 \gamma_6^2 + 16) \gamma_6^4,
	
	\\
	\noalign{\medskip}
	
	\frac{24\gamma_6^6}{\lambda^2+\gamma_6^2}\, \bach_{44}=  
	-4 ( \gamma_6^2 - 4) \lambda^4
	-(4 \gamma_6^4 + 5 \gamma_6^2 - 48 ) \gamma_6^2 \lambda^2
	+(12 \gamma_6^4 - 33 \gamma_6^2 + 48) \gamma_6^4,
\end{array}
\]
which implies 
\[  
\begin{array}{l}
	\frac{3\gamma_6^6}{\lambda^2(\lambda^2+\gamma_6^2)} (3\,\bach_{11} +\bach_{44}) =  
	(7 \gamma_6^2 + 8) \lambda^2 + 
	4  (\gamma_6^4 - \gamma_6^2 + 6) \gamma_6^2 .
\end{array}
\]
Since $\lambda\gamma_6\neq 0$,  we conclude   that the Bach tensor does not vanish.

%
%

\subsection{The structure operator is diagonalizable of rank one with Lorentzian kernel}\label{se:Lorentz-2}

In this setting, it is possible to choose    an orthonormal basis $\{e_1,e_2,e_3,e_4\}$ of $\mathfrak{g}=\mathfrak{h}_3\rtimes \mathfrak{r}$, with $e_3$ timelike,  where  $\mathfrak{h}_3=\spn\{e_1,e_2,e_3\}$ and $\mathfrak{r}=\spn\{e_4\}$, so that 
the left-invariant metrics are described by
\[
\begin{array}{ll}
	[e_1,e_3]=-\lambda e_2,  &
	[e_1,e_4]=\gamma_1 e_1+\gamma_2 e_2 + \gamma_3 e_3,  
	\\
	\noalign{\medskip}
	[e_2,e_4]=\gamma_4 e_2,  &
	[e_3,e_4]=\gamma_5 e_1+\gamma_6 e_2-(\gamma_1 -\gamma_4)e_3,
\end{array}
\] 
where $\lambda\neq 0$ and $\gamma_1$, $\dots$, $\gamma_6\in\mathbb{R}$. 
In this case  we study the vanishing of the  Bach tensor   analyzing      the cases $\gamma_4=0$ and $\gamma_4\neq 0$ separately.

\smallskip
 
\subsubsection{Case $\gamma_4= 0$}\label{se:Lorentz-2-1}

Since $\lambda\neq 0$, we consider the orthogonal basis $\hat e_i=\frac{1}{\lambda}e_i$ so that we can assume $\lambda=1$ working in the homothetic class of the initial metric.  In this case, we compute a  Gr\"obner basis  of the ideal  $\langle \bach_{ij}\rangle\subset \mathbb{R}[\gamma_2$, $\gamma_3$, $\gamma_5$, $\gamma_6$, $\gamma_1]$  with respect to the   lexicographical order (consisting of 25 polynomials) and get  that the polynomials
\[
 \mathbf{g}_{1}=  \gamma_1^4,\quad
 \mathbf{g}_{2}= \gamma_3^4
 \quad\text{and}\quad
 \mathbf{g}_{3}= \gamma_5^4
\]
belong to the ideal. Thus, necessarily $\gamma_1=\gamma_3=\gamma_5=0$  and a direct calculation shows   
$\bach_{22} = -\frac{5}{6} (\gamma_6^2-\gamma_2^2+1)^2$, which implies $\gamma_2 = \varepsilon \sqrt{\gamma_6^2+1}$, where $\varepsilon^2=1$. Under this assumption the Bach tensor vanishes and  the corresponding left-invariant metric is given by
\[ 
	[e_1,e_3]=-  e_2,  \qquad 
	[e_1,e_4]=   \varepsilon \sqrt{\gamma_6^2+1} \, e_2  ,   \qquad  
	[e_3,e_4]= \gamma_6 e_2 . 
\]
Moreover, the Ricci operator is $2$-step nilpotent and   
$\xi = -\gamma_6  e_1 + \varepsilon \sqrt{\gamma_6^2+1}\,   e_3 +   e_4$ 
is a null parallel vector field   so that     $R(x,y)=0$ and $\nabla_xR=0$ for all $x,y\in \xi^\perp=\spn\{\xi$, $ e_2$, $e_1+\gamma_6 e_4\}$. Thus, the underlying structure is a plane wave.
Considering the basis 
$$
\te_1=e_1,\quad \te_2=e_2,\quad \te_3=-e_3,\quad 
\te_4=-\gamma_6 e_1+\varepsilon\sqrt{\gamma_6^2+1} \, e_3+e_4,
$$
one has that the only non-zero bracket is $[\te_1,\te_3]=\te_2$. Hence the underlying group structure is the product $H_3\times\mathbb{R}$ and the metric is isomorphically homothetic to the Bach-flat metric discussed in Remark~\ref{re:producto}.

%

\smallskip
 
\subsubsection{Case $\gamma_4\neq 0$}\label{se:Lorentz-2-2}
 
Taking the orthogonal basis $\hat e_i=\frac{1}{\gamma_4} e_i$ we can assume $\gamma_4=1$ working in the homothetic class of the initial metric. We consider the ideal $\langle \bach_{ij}\rangle$ in the polynomial  ring $\mathbb{R}[\gamma_2$, $\gamma_3$, $\gamma_5$, $\gamma_6$, $\gamma_1$, $\lambda]$ and use the lexicographical order to compute a   Gr\"obner basis. As a consequence, we get 61 polynomials, being two of them
\[
\begin{array}{l}
	\mathbf{g}_1 =  
	 (\gamma_3 + \gamma_5) \gamma_6 (\lambda^2 + 1) (\lambda^2 + 
	 25) (4 \lambda^2 + 1) 
	 \quad\text{and}\quad
	 \\
	 \noalign{\medskip}
	 \mathbf{g}_2 =
	 (\gamma_3 + \gamma_5) (7 \lambda^2 - 2)(\lambda^2 + 1)
	  (\lambda^2 + 2) (\lambda^2 + 25) (4 \lambda^2 + 1) .
\end{array}
\]
Hence, either $\gamma_3+\gamma_5=0$ or $\gamma_3+\gamma_5\neq0$, $\gamma_6=0$, $7 \lambda^2 - 2=0$. Next we analyze these two cases by separate.

\smallskip
 
\paragraph{\underline{\emph{Case $\gamma_3+\gamma_5=0$}}}\label{se:Lorentz-2-2-1}

Continuing  in the polynomial  ring $\mathbb{R}[\gamma_2$, $\gamma_3$, $\gamma_5$, $\gamma_6$, $\gamma_1$, $\lambda]$ with  the lexicographical order, for the ideal $\langle \bach_{ij}\cup\{\gamma_3+\gamma_5\}\rangle$ we obtain the polynomials
\[
\begin{array}{l}
	\mathbf{g}'_1 =  
	-((\gamma_1 - 1) \lambda^2 - \gamma_6^2) (\lambda^2 + 1) 
	(2 \lambda^2 + 1) (4 \lambda^2 + 1) (5 \lambda^2 + 8)
	 (25 \lambda^2 + 16) 
	\quad\text{and}\quad
	\\
	\noalign{\medskip}
	\mathbf{g}'_2 =
	((\gamma_1 - 1) \gamma_2 - \gamma_5 \gamma_6) (2 \lambda^2 + 
	1) (5 \lambda^2 + 8) (25 \lambda^2 + 16).
\end{array}
\]	
Hence, $\mathbf{\tilde g}'_1 = (  \gamma_1-1) \lambda^2 - \gamma_6^2$ and $\mathbf{\tilde g}'_2 = (\gamma_1 - 1) \gamma_2 - \gamma_5 \gamma_6$ must vanish. Repeating the process for the ideal $\langle \bach_{ij}\cup\{\gamma_3+\gamma_5$, $\mathbf{\tilde g}'_1$, $\mathbf{\tilde g}'_2\}\rangle$ we get the polynomials
\[
\begin{array}{l}
	\mathbf{g}'_3 = 
		( \gamma_5^2 -  \gamma_1 ( \gamma_1 - 1)) (5 \lambda^2 + 8),
		\\
		\noalign{\medskip}
	\mathbf{g}'_4 = 
		8 ( \gamma_5^2 -  \gamma_1 ( \gamma_1 -  
		1))  +  \gamma_2^2-  \gamma_1 \lambda^2 
		\quad\text{and}\quad
		\\
		\noalign{\medskip}
	\mathbf{g}'_5 = 
			8  \gamma_5 (  \gamma_5^2 -  \gamma_1 ( \gamma_1 - 
			1) ) +  \gamma_2  \gamma_6 -  \gamma_5  \lambda^2 .
\end{array}
\]
Thus, $  \gamma_2^2-  \gamma_1 \lambda^2 =0$ and $  \gamma_2  \gamma_6-  \gamma_5  \lambda^2 =0$, and these relations, together with $\gamma_3+\gamma_5=0$ and $(  \gamma_1-1) \lambda^2 - \gamma_6^2=0$, imply
\[
	\gamma_3=-\gamma_5,\quad
	\gamma_1 = \tfrac{\lambda^2+\gamma_6^2}{\lambda^2},\quad
	\gamma_2=\varepsilon \sqrt{\lambda^2+\gamma_6^2},\quad
	\gamma_5= \tfrac{\gamma_2\gamma_6}{\lambda^2},
\] 
where $\varepsilon^2=1$. A straightforward calculation shows that the Bach tensor vanishes for the corresponding left-invariant metric, given by
\[
\begin{array}{ll}
	[e_1,e_3]=-\lambda e_2,  &
	[e_1,e_4]=\frac{\lambda^2+\gamma_6^2}{\lambda^2} e_1
				+\varepsilon \sqrt{\lambda^2+\gamma_6^2} \,e_2 
				-\frac{\varepsilon \gamma_6\sqrt{\lambda^2+\gamma_6^2}}{\lambda^2} e_3,  
	\\
	\noalign{\medskip}
	[e_2,e_4]=  e_2,  &
	[e_3,e_4]=\frac{\varepsilon \gamma_6\sqrt{\lambda^2+\gamma_6^2}}{\lambda^2} e_1
	 		+\gamma_6 e_2- \frac{ \gamma_6^2}{\lambda^2}   e_3.
\end{array}
\] 
Now, considering a new basis defined by
\[
	\te_1 = \tfrac{\varepsilon\sqrt{\lambda^2+\gamma_6^2}}{\lambda}\,e_1 -\tfrac{\gamma_6}{\lambda}e_3,\quad
	\te_2 = e_2,\quad
	\te_3 = -\tfrac{\gamma_6}{\lambda}e_1
	 + \tfrac{\varepsilon\sqrt{\lambda^2+\gamma_6^2}}{\lambda}\,e_3,\quad
	 \te_4 = e_4,
\]
the inner product remains invariant while the Lie bracket transforms into
\[
	[\te_1,\te_3] = -\lambda \te_2,\qquad
	[\te_1,\te_4] = \te_1 + \lambda \te_2,\qquad
	[\te_2,\te_4] = \te_2,
\]
and a straightforward calculation shows that this metric is not locally symmetric, locally conformally flat or Einstein.
Moreover, $(\te_1,\te_2,\te_3,\te_4)\mapsto(\te_1,-\te_2,\te_3,\te_4)$ defines an isometry interchanging $\lambda$ and $-\lambda$, and therefore we can restrict the parameter $\lambda$ to $\lambda>0$. This case corresponds to the metrics in Theorem~\ref{th:ce}--(L.i).

\smallskip
 
\paragraph{\underline{\emph{Case $\gamma_3+\gamma_5\neq0$, $\gamma_6=0$, $7\lambda^2-2=0$}}}\label{se:Lorentz-2-2-2}
 
For this setting we consider  the ideal $\langle \bach_{ij}  \cup \{7\lambda^2-2\}\rangle$ in the polynomial ring
$\mathbb{R}[\gamma_2$, $\gamma_3$, $\gamma_5$, $\gamma_1$, $\lambda]$. In this case, we use the graded reverse lexicographical order to see that the polynomials 
\[
\begin{array}{l}
	\mathbf{g}''_1 = \gamma_2 \gamma_5,
	\\
	\noalign{\medskip}
	\mathbf{g}''_2 = \gamma_2 \gamma_3, 
	\\
	\noalign{\medskip}
	\mathbf{g}''_3 = 8  \gamma_1^2 + 21  \gamma_2^2 - 8  \gamma_1 + 
	8  \gamma_3  \gamma_5 - 6
	\quad\text{and}\quad
	\\
	\noalign{\medskip}
	\mathbf{g}''_4 = -22  \gamma_1^2 + 21  \gamma_3^2 + 21  \gamma_5^2 + 22  \gamma_1 + 
	20  \gamma_3  \gamma_5
\end{array}
\]
belong to the ideal. Since $\gamma_3+\gamma_5\neq 0$, clearly $\gamma_2=0$. Now, 
$\mathbf{g}''_3 =0$ and $\mathbf{g}''_4=0$, together with  $7\lambda^2-2=0$, imply
\[
	\lambda = \varepsilon_1 \sqrt{\tfrac{2}{7}},\quad
	\gamma_5 = \varepsilon_2 \sqrt{\tfrac{11}{14}}-\gamma_3,\quad
	\gamma_1 = \tfrac{1}{2}+\varepsilon_3 \sqrt{\gamma_3^2+1-\varepsilon_2\sqrt{\tfrac{11}{14}}\,\gamma_3},
\]
where $\varepsilon_1^2=\varepsilon_2^2=\varepsilon_3^2=1$. A straightforward calculation shows that the corresponding left-invariant metric, given by 
\[
\begin{array}{l}
	[e_1,e_3]=-\varepsilon_1 \sqrt{\frac{2}{7}}\, e_2,  \qquad
	[e_1,e_4]=
		\left(\frac{1}{2}+\varepsilon_3 \sqrt{\gamma_3^2+1-\varepsilon_2\sqrt{\frac{11}{14}}\,\gamma_3}\right)  e_1  
		+ \gamma_3 e_3,  
	\\
	\noalign{\medskip}
	[e_2,e_4]=  e_2, \quad\quad 
	[e_3,e_4]= \left(\varepsilon_2 \sqrt{\frac{11}{14}}-\gamma_3\right) e_1 
				+ \left(  
					\frac{1}{2}
					-\varepsilon_3 \sqrt{\gamma_3^2+1-\varepsilon_2\sqrt{\frac{11}{14}}\,\gamma_3}
					\right)e_3,
\end{array}
\] 
is Bach-flat. 

Let $\gamma_3'=\sqrt{\frac{14}{11}}\varepsilon_2\gamma_3$.
A crucial observation is that it is possible to reduce this case eliminating the parameter $\gamma_3'$. We proceed as follows.  First, we consider a new basis $\{\te_i\}$ defined by
\[
\te_1 = 2\sqrt{14}\,e_1,\qquad
\te_2 = 2\sqrt{14}\,\varepsilon_1\varepsilon_2 e_2,\qquad
\te_3 = 2\sqrt{14}\,\varepsilon_2 e_3,\qquad
\te_4 =  2\sqrt{14}\,e_4.
\]
Then, the Lie bracket transforms into
\[
\begin{array}{l}
	[\te_1,\te_3]=-4 \te_2,  \qquad
	[\te_2,\te_4]=  2\sqrt{14}\,\te_2, 
	
	\\
	\noalign{\medskip}
		
	[\te_1,\te_4]=
		\left(\sqrt{14}
		+2\varepsilon_3\sqrt{
		 11 \gamma_3'(\gamma_3'- 1)
		+14}
		\right)\te_1
		+2\sqrt{11}\gamma_3' \te_3 ,	 
	
	\\
	\noalign{\medskip}
		
	[\te_3,\te_4]= 
		-2\sqrt{11}(\gamma_3'-1)\te_1
		+\left(\sqrt{14}
				-2\varepsilon_3\sqrt{
				 11 \gamma_3'(\gamma_3'- 1)
				+14}
				\right)\te_3 ,
\end{array}
\]
while the inner product is given by $56\langle\cdot,\cdot\rangle$. Since we are working at the homothetic level, we can maintain the
initial inner product $\langle\cdot,\cdot\rangle$  remaining in the same homothetic class. 

Secondly,  we make another change of basis taking
$\{f_i\}$ given by
\[
\f_1=\te_4,\qquad
\f_2=\tfrac{\tf_2}{\|\tf_2\|},\qquad
\f_3=\tfrac{\tf_3}{\|\tf_3\|},\qquad
\f_4= \te_2,
\]
where
\[
\tf_2 = 
\left\{
\begin{array}{ll}
	-2\sqrt{\frac{14}{11}}\,\te_1+\te_3
	& \text{if $(\varepsilon_3,\gamma_3')=(1,0)$},
	
	\\
	\noalign{\medskip}
	
	\te_1
	& \text{if $(\varepsilon_3,\gamma_3')=(1,1)$},
	
	\\
	\noalign{\medskip}
	
	\te_3
		& \text{if $(\varepsilon_3,\gamma_3')=(-1,0)$},
		
	\\
	\noalign{\medskip}

	-\frac{1}{2}\sqrt{\frac{11}{14}}\,\te_1+\te_3
	& \text{if $(\varepsilon_3,\gamma_3')=(-1,1)$},
	
	\\
	\noalign{\medskip}
		
	\tfrac{\sqrt{14}
				+\varepsilon_3\sqrt{
				 11 \gamma_3'(\gamma_3' -1)
				+14}}{\sqrt{11}(\gamma_3'-1)}\, \te_1 + \te_3
	& \text{if $\gamma_3'\notin\{0,1\}$,}
\end{array}
\right.
\]
and
\[
\tf_3 = 
\left\{
\begin{array}{ll}
	-\frac{1}{2}\sqrt{\frac{11}{14}}\,\te_1+\te_3
	& \text{if $(\varepsilon_3,\gamma_3')=(1,0)$},
	
	\\
	\noalign{\medskip}
	
	\te_3
	& \text{if $(\varepsilon_3,\gamma_3')=(1,1)$},
	
	\\
	\noalign{\medskip}
	
	\te_1
		& \text{if $(\varepsilon_3,\gamma_3')=(-1,0)$},
		
	\\
	\noalign{\medskip}

	-2\sqrt{\frac{14}{11}}\,\te_1+\te_3
	& \text{if $(\varepsilon_3,\gamma_3')=(-1,1)$},
	
	\\
	\noalign{\medskip}
		
	-\tfrac{\sqrt{14}
					-\varepsilon_3\sqrt{
					 11 \gamma_3'(\gamma_3' -1)
					+14}}{\sqrt{11}\,\gamma_3'}\, \te_1 + \te_3
	& \text{if $\gamma_3'\notin\{0,1\}$.}
\end{array}
\right.
\]
A long but straightforward calculation shows that
the inner product transforms into $\operatorname{diag}[1,\varepsilon_3,-\varepsilon_3,1]$, while  the Lie bracket is given by
\[
\begin{array}{ll}
[f_1,f_2] = -3 \sqrt{14} f_2 + 2 \delta_1  \sqrt{11} f_3,
&
\qquad
[f_1,f_3] = \sqrt{14}f_3,

\\
\noalign{\medskip}

[f_1,f_4] = -2\sqrt{14} f_4,
&
\qquad
[f_2,f_3] = 4\delta_2 f_4,
\end{array}
\]
where  
\[
(\delta_1,\delta_2)=\left\{
\begin{array}{ll}
	(1,1)
	& \text{if $(\varepsilon_3,\gamma_3')=(1,0)$},
	
	\\
	\noalign{\medskip}
	
	(-1,-1)
	& \text{if $(\varepsilon_3,\gamma_3')=(1,1)$},
	
	\\
	\noalign{\medskip}
	
	(-1,1)
		& \text{if $(\varepsilon_3,\gamma_3')=(-1,0)$},
		
	\\
	\noalign{\medskip}

	(1,-1)
	& \text{if $(\varepsilon_3,\gamma_3')=(-1,1)$},
	
	\\
	\noalign{\medskip}
		
	(\operatorname{sign}(1-\gamma_3'),
	\operatorname{sign}(1-\gamma_3'))
	& \text{if $\varepsilon_3=1$ and $\gamma_3'\notin\{0,1\}$},

	\\
	\noalign{\medskip}
			
	(\operatorname{sign}(\gamma_3'),
	-\operatorname{sign}(\gamma_3'))
	& \text{if $\varepsilon_3=-1$ and $\gamma_3'\notin\{0,1\}$},
\end{array}
\right.
\]
so we have $\delta_1^2=\delta_2^2=1$.

Now, a final calculation shows that the metric above is  not locally symmetric, locally conformally flat or Einstein. 
One may further assume $\delta_2=~1$, just considering the isometry given by $(f_1,f_2,f_3,f_4)\mapsto(f_1,f_2,\delta_2 f_3,f_4)$. This corresponds to the metric Theorem~\ref{th:strictly Bach flat}--(i).

 \smallskip

\subsection{The structure operator is $2$-step nilpotent}\label{se:Lorentz-3} 

In this last case,   there exists a pseudo-orthonormal basis $\{u_1,u_2,u_3,u_4\}$ of $\mathfrak{g}=\mathfrak{h}_3\rtimes \mathfrak{r}$, with $\langle u_1,u_2\rangle=\langle u_3,u_3\rangle=\langle u_4,u_4\rangle=1$,  where  $\mathfrak{h}_3=\spn\{u_1,u_2,u_3\}$ and $\mathfrak{r}=\spn\{u_4\}$, so that 
\[
\begin{array}{ll}
	[u_1,u_3]=-\varepsilon u_2,&
	[u_1,u_4]=\gamma_1 u_1+\gamma_2 u_2+\gamma_3 u_3,   
	
	\\
	\noalign{\medskip}
	
	[u_2,u_4]=\gamma_4 u_2, &
	[u_3,u_4]=\gamma_5 u_1+\gamma_6 u_2 - (\gamma_1-\gamma_4)u_3 ,
\end{array}
\] 
with $\varepsilon^2=1$ and $\gamma_1$, $\dots$, $\gamma_6\in\mathbb{R}$.  
A direct calculation shows that $\bach_{34} = -\frac{3}{4}\varepsilon \gamma_4\gamma_5^2$. We analyze by separate the cases $\gamma_4=0$ and $\gamma_4\neq0$, $\gamma_5=0$.

\smallskip

\subsubsection{Case $\gamma_4=0$}\label{se:Lorentz-3-1}

In this case, we have
\[
\bach_{44}= \tfrac{1}{8} (3\gamma_1^2 + (\gamma_3+4\gamma_6)\gamma_5)^2 .
\]

If $\gamma_5=0$ then necessarily $\gamma_1=0$ and a straightforward calculation shows that the Bach-flatness condition is satisfied. Hence,  the corresponding left-invariant metric is given by
\[ 
	[u_1,u_3]=-\varepsilon u_2,\qquad
	[u_1,u_4]= \gamma_2 u_2+\gamma_3 u_3,\qquad
	[u_3,u_4]= \gamma_6 u_2,
\]
and  $u_2$ is a null parallel vector field so that   the curvature tensor satisfies $R(x,y)=0$ and $\nabla_xR=0$ for all $x,y\in u_2^\perp = \langle u_2$, $u_3$, $u_4\rangle$. Moreover, the only non-zero component of the Ricci tensor is $\rho_{11}=-\frac{1}{2}(\gamma_3^2-\gamma_6^2)$ which implies that the Ricci operator is isotropic  and therefore the underlying structure is a plane wave. Furthermore the metric is locally conformally flat or locally symmetric if and only if $\gamma_3(\gamma_3+\gamma_6)=0$. This corresponds to Theorem~\ref{th:plane-waves}-(L).


\medskip

If $\gamma_5\neq 0$, then $\gamma_6= -\frac{3\gamma_1^2 +\gamma_3\gamma_5}{4\gamma_5}$ and the only non-zero component of the Bach tensor is determined by
\[
\begin{array}{l}
	\frac{64\gamma_5^2}{3}\bach_{11} = 
	16  \gamma_5^4 \gamma_2^2  -  
	8 \gamma_1 (3 (\gamma_1^2 + \gamma_3 \gamma_5) - 2 \gamma_1^2) \gamma_5^2  \gamma_2   
	
	\\
	\noalign{\medskip}
	\phantom{\frac{64\gamma_5^2}{3}\bach_{11}=}
	
	+  (\gamma_1^2 + 
	\gamma_3 \gamma_5) (8 (\gamma_1^2 + \gamma_3 \gamma_5)^2 + 
	9 (\gamma_1^2 + \gamma_3 \gamma_5) \gamma_1^2 - 
	12 \gamma_1^4) + 4 (\gamma_1^6 + 4 \gamma_5^4),
\end{array}
\]
which is a polynomial of  degree two in the variable $\gamma_2$ whose discriminant reduces to  $\Delta= -512 \gamma_5^4 (2\gamma_5^4 + (\gamma_1^2 + \gamma_3 \gamma_5)^3)$. Hence, 
if $\gamma_1^2 + \gamma_3 \gamma_5\geq 0$ we have  $\Delta <0 $ and the metric cannot be Bach-flat.
Thus, necessarily $\gamma_1^2 + \gamma_3 \gamma_5<0$ and, in particular, $\gamma_3 \gamma_5<0$.

Let $\gamma_1'$ and $\gamma_5'$  be new  variables defined by
$\gamma_1'=\frac{\gamma_1}{\gamma_3^2}$ and 
$\gamma_5'=\frac{\gamma_5}{\gamma_3^3}$ so that the condition $\gamma_1^2 + \gamma_3 \gamma_5<0$ reads  $\gamma_3^4((\gamma_1')^2+\gamma_5')<0$. In order to simplify the writing, let $\kappa=-((\gamma_1')^2 +  \gamma_5')$, $\kappa>0$. We make a change of basis and from now on we consider   a new basis $\{ \tu_i \}$ defined by
\[
	\tu_1= \tfrac{\varepsilon}{\gamma_3\sqrt{\kappa}} u_1,
	\qquad
	\tu_2=  \tfrac{\varepsilon}{\gamma_3^3\sqrt{\kappa}} u_2,
	\qquad
	\tu_3= \tfrac{\varepsilon}{\gamma_3^2\sqrt{\kappa}} u_3,
	\quad
	\tu_4= \tfrac{1}{\gamma_3^2\sqrt{\kappa}} u_4.
\]
Hence, the Lie bracket transforms into 
%
%
\begin{equation}\label{eq:Lie bracket}
\begin{array}{l}
	[\tu_1,\tu_3] = 			
		-\tfrac{1}{\sqrt{\kappa}}\tu_2,
	\qquad
	[\tu_1,\tu_4] = 			
		\tfrac{1}{\sqrt{\kappa}}
		(\gamma_1' \tu_1+\gamma_2\tu_2+\tu_3),
		
	\\
	\noalign{\medskip}
	
	[\tu_3,\tu_4] = 			
			\tfrac{1}{\sqrt{\kappa}}
			\left(
			\gamma_5' \tu_1
			-\frac{3(\gamma_1')^2+\gamma_5'}{4\gamma_5'}\tu_2
			-\gamma_1'\tu_3
			\right),
	
\end{array}
\end{equation}
while the inner product is given by $\frac{1}{\gamma_3^4\kappa}\langle\cdot,\cdot\rangle$. Since we are working at the homothetic level  we can maintain the initial inner product   remaining in the same homothetic class. 
With respect to the new basis  the Bach tensor is determined by
\[
\begin{array}{l}
	\frac{64(\gamma_5')^2\kappa^2}{3}\bach_{11} = 16  (\gamma_5')^4 \gamma_2^2  
	+  8 \gamma_1' (\gamma_5')^2 
	\left( 2 (\gamma_1')^2+3 \kappa \right) \gamma_2   
	
	\\
	\noalign{\medskip}
	\phantom{\frac{64(\gamma_5')^2\kappa^2}{3} \bach_{11}=}
	
	+  
	(12 (\gamma_1')^4+ 9\kappa (\gamma_1')^2 -8 \kappa^2)\kappa 
	+ 4 ((\gamma_1')^6 + 4 (\gamma_5')^4)
\end{array}
\]
and analyzing the Ricci operator of the metric, given by
\[ 
	\Ricci=\left(
	\begin{array}{cccc}
		\frac{2\gamma_5'+3\kappa}{8\kappa}&
		\frac{(\gamma_5')^2}{2\kappa}&
		-\frac{\gamma_1'\gamma_5'}{2\kappa}&
		0
		
		\\
		\noalign{\medskip}

		-\frac{4(8\gamma_1'\gamma_2+3)(\gamma_5')^2
		-12\gamma_5'\kappa
		-9\kappa^2}{32(\gamma_5')^2\kappa}
		&
		\frac{2\gamma_5'+3\kappa}{8\kappa} &
		-\frac{2\gamma_2(\gamma_5')^2+3\gamma_1'\kappa}{4\gamma_5'\kappa}   &
		-\frac{\gamma_5'}{2\kappa}
		
		\\
		\noalign{\medskip}
		
		-\frac{2\gamma_2(\gamma_5')^2+3\gamma_1'\kappa}{4\gamma_5'\kappa}   &
		-\frac{\gamma_1'\gamma_5'}{2\kappa} &
		-\frac{ 2\gamma_5'+3\kappa }{4\kappa} &
		0
		
		\\
		\noalign{\medskip}
		
		-\frac{\gamma_5'}{2\kappa} &
		0 &
		0 &
		\frac{3}{4}
		
	\end{array}
	\right),
	\]
it follows that 
	$
	\bach_{11}=
	\frac{8\kappa}{(\gamma_5')^2}\det(\Ricci)+\frac{(\gamma_5')^2}{\kappa^2}
	$.
	Hence, if $\gamma_5'\neq 0$, $\kappa=-((\gamma_1')^2 +  \gamma_5')>0$ and $\det(\Ricci)=-\frac{(\gamma_5')^4}{8\kappa^3}$, which is equivalent to
	\[
		\gamma_2 =\tfrac{1}{4(\gamma_5')^2} \left(
		\gamma_1'(2\gamma_5'- \kappa)
		+ 2 \varepsilon' \sqrt{2 \left( \kappa^3-2(\gamma_5')^4 \right)} 
		\right)
		\,\,\text{with}\,\,
		\kappa^3-2(\gamma_5')^4\geq0 ,\,
		(\varepsilon')^2=1,
	\]
then the left-invariant  metric given by Equation~\eqref{eq:Lie bracket} is Bach-flat and a straightforward calculation shows that it has positive scalar curvature $\tau=\frac{3}{4}$ and it  is not locally symmetric, locally conformally flat or Einstein.  
	Note that the conditions $\gamma_5'\neq 0$, $\kappa=-((\gamma_1')^2 +  \gamma_5')>0$ and $\kappa^3-2(\gamma_5')^4\geq0$ are equivalent to  $-\frac{1}{2}\leq\gamma_5'<0$,
	$(\gamma_1')^2\leq -\gamma_5'\left(
	\sqrt[3]{2\gamma_5'}+1\right)$.
	These metrics correspond to those in Theorem~\ref{th:strictly Bach flat}--(ii).

\smallskip

\subsubsection{Case $\gamma_4\neq0$, $\gamma_5=0$}\label{se:Lorentz-3-2}

In this case, $\bach_{44}=\frac{1}{24}(3\gamma_1-\gamma_4)^2
(\gamma_1-\gamma_4)(3\gamma_1+\gamma_4)
$. Hence, $\gamma_4=3\gamma_1$, $\gamma_4=\gamma_1$ or $\gamma_4=-3\gamma_1$. 
If $\gamma_4=3\gamma_1$ then a direct calculation shows that the Bach-flatness condition is equivalent to $\gamma_6=-\gamma_3$ and the corresponding metric is locally conformally flat. Moreover, the metric is both locally symmetric and Einstein if and only if $\gamma_2=0$.
Next we analyze by separate the other two cases,  $\gamma_4=\gamma_1$ and  $\gamma_4=-3\gamma_1$.

\medskip

\paragraph{\underline{\emph{Case $\gamma_4=\gamma_1$}}}\label{se:Lorentz-3-2-1}

A direct calculation shows that the non-zero components of the Bach tensor are determined by
\[
	\bach_{11}=\tfrac{1}{4}\gamma_1^2(\gamma_3-7\gamma_6)(\gamma_3-\gamma_6)
	\quad\text{and}\quad
	\bach_{13}=\tfrac{3}{4}\gamma_1^3 (\gamma_3-\gamma_6).
\]
Since $\gamma_1\neq 0$, it follows that $\gamma_6=\gamma_3$ and the left-invariant Bach-flat metric is given by
\[
\begin{array}{ll}
	[u_1,u_3]=-\varepsilon u_2,&
	[u_1,u_4]=\gamma_1 u_1+\gamma_2 u_2+\gamma_3 u_3,   
	
	\\
	\noalign{\medskip}
	
	[u_2,u_4]=\gamma_1 u_2, &
	[u_3,u_4]=\gamma_3 u_2.
\end{array}
\]
Let $\gamma_2' = \frac{\gamma_1\gamma_2+\gamma_3^2}{|\gamma_1|}\in\mathbb{R}$ and let $\nu=-1,1, $ depending on whether $\gamma_1$ is negative or positive. Define a new basis $\{\tu_i\}$ by
\[ 
\begin{array}{ll}
	\tu_1
	= -\tfrac{2}{|\gamma_1|^\frac{1}{2}}u_1 + \tfrac{\gamma_3^2}{|\gamma_1|^\frac{5}{2}}u_2
	-\tfrac{2\nu\gamma_3}{|\gamma_1|^\frac{3}{2}}u_3  , & \quad  
	\tu_2 =   -\tfrac{8}{|\gamma_1|^\frac{3}{2}} u_2,
	\\
	\noalign{\medskip}
	\tu_3=
	-\tfrac{4\nu\varepsilon\gamma_3}{\gamma_1^2} u_2 + \tfrac{4\varepsilon}{|\gamma_1|}u_3,  &
	\quad
	\tu_4=   \tfrac{4}{\gamma_1} u_4.	 
\end{array}
\] 
Now, a direct calculation shows that  the Lie bracket transforms into
\[
[\tu_1,\tu_3]= -\tu_2,\qquad
[\tu_1,\tu_4]=  4 \tu_1 + \gamma_2' \tu_2, \qquad
[\tu_2,\tu_4]=  4  \tu_2 ,
\]
while  the inner product     is given by $\frac{16}{\gamma_1^2} \langle\cdot,\cdot\rangle$. Since we are working at the homothetic level we can maintain the initial inner product   remaining in the same homothetic class. A straightforward calculation shows that the above metrics are not locally symmetric, locally conformally flat or Einstein. They correspond to the metrics in Theorem~\ref{th:ce}--(L.ii).

\medskip

\paragraph{\underline{\emph{Case $\gamma_4=-3\gamma_1$}}}\label{se:Lorentz-3-2-1}

The Bach tensor is determined by
\[
	\bach_{11}=-\tfrac{1}{4}\gamma_1^2
	(7\gamma_3^2+41\gamma_6^2+32\gamma_1\gamma_2-24\gamma_3\gamma_6)\,,
	\quad\quad
	\bach_{13}=\tfrac{7}{4}\gamma_1^3 (\gamma_3-5\gamma_6),
\]
and since $\gamma_1\neq 0$ we get $\gamma_3=5\gamma_6$ and
$\gamma_2=-\frac{3\gamma_6^2}{\gamma_1}$.
Thus, the Bach-flat left-invariant metric corresponds to
\[
\begin{array}{ll}
	[u_1,u_3]=-\varepsilon u_2,&
	[u_1,u_4]=\gamma_1 u_1-\frac{3\gamma_6^2}{\gamma_1} u_2+5\gamma_6 u_3,   
	
	\\
	\noalign{\medskip}
	
	[u_2,u_4]=-3\gamma_1 u_2, &
	[u_3,u_4]=\gamma_6 u_2 - 4\gamma_1 u_3 .
\end{array}
\]
Let $\nu=-1$,   $1$, depending on whether $\gamma_1$ is negative  or positive. Considering the basis  $\{\tu_i\}$ defined by 
\[ 
\begin{array}{ll}
\tu_1
= -\tfrac{1}{|\gamma_1|^\frac{3}{2}}\left(\gamma_1 u_1 -\tfrac{\gamma_6^2}{2\gamma_1}u_2+\gamma_6 u_3\right) , & \quad  
\tu_2 =   -\tfrac{1}{\gamma_1 |\gamma_1|^\frac{1}{2}} u_2,
\\
\noalign{\medskip}
\tu_3=
-\tfrac{\nu\varepsilon\gamma_6}{\gamma_1^2} u_2 + \tfrac{\varepsilon}{|\gamma_1|}u_3,  &
\quad
\tu_4=   -\tfrac{\nu}{|\gamma_1|} u_4 ,	 
\end{array}
\] 
a direct calculation shows that  the Lie bracket transforms into
\[
[\tu_1,\tu_3]= -\tu_2,\qquad
[\tu_1,\tu_4]=  - \tu_1, \qquad
[\tu_2,\tu_4]=  3  \tu_2,\qquad
[\tu_3,\tu_4]= 4 \tu_3,
\]
while  the inner product     is given by $\frac{1}{\gamma_1^2} \langle\cdot,\cdot\rangle$. Since we are working at the homothetic level  we can maintain the initial inner product   remaining in the same homothetic class. Now, a straightforward calculation shows that the above metric is not locally symmetric, locally conformally flat or Einstein. It corresponds to the metric in Theorem~\ref{th:ce}--(L.iii).

\section{Proof of Theorem~\ref{th:ce}}\label{se:proof-1} 

We consider the Bach-flat metrics obtained in Sections $\S$\ref{se:2}, $\S$\ref{se:Riemann}  and $\S$\ref{se:EGR-Lorentzian}, except those which are $pp$-waves already covered by Remark~\ref{re:producto} and Theorem~\ref{th:plane-waves}. Next we show that the metric is conformally Einstein in all cases of Theorem~\ref{th:ce} by analyzing each situation separately.

Recall that a metric is conformally Einstein if and only if there exists a (locally defined) nowhere zero function $\varphi$ so that $\overline{g}=\varphi^{-2}g$ satisfies Equation \eqref{eq:hoy2}. Moreover, setting $\varphi=e^{\sigma}$ one has that the gradient of the function $\sigma$ satisfies the conformal Cotton-flat equation 
\eqref{eq:hoy3} since $\overline{g}=e^{-2\sigma}g$ is Cotton-flat. 
Set $\xi=\nabla\sigma$.
Now a straightforward calculation shows that $\nabla\varphi=\varphi\xi$ and
$$
\operatorname{Hes}_\varphi(X,Y)=\varphi\{ \langle X,\xi\rangle\langle Y,\xi\rangle+\langle\nabla_X\xi,Y\rangle \}.
$$
In order to analyze Equation \eqref{eq:hoy2} we consider the symmetric $(0,2)$-tensor field
$$
\begin{array}{rcl}
\mathfrak{CE}(X,Y)&=&2\operatorname{Hes}_\varphi(X,Y)+\varphi\rho(X,Y)-\frac{1}{4}\{ 2\Delta\varphi+\varphi\tau\}\langle X,Y\rangle
\\
\noalign{\medskip}
&=&
2\varphi\{ \langle X,\xi\rangle\langle Y,\xi\rangle+\langle\nabla_X\xi,Y\rangle \}+\varphi\rho(X,Y) -
\frac{1}{4}\{2\Delta\varphi+\varphi\tau\} \langle X,Y\rangle 
\end{array}
$$
and evaluate it on the left-invariant vector fields obtained from the corresponding basis of each Lie algebra in Theorem~\ref{th:ce}.

\subsubsection*{Proof of Theorem~\ref{th:ce}-(D.i)}
Let $\{U_1,U_2,U_3,U_4\}$ be the pseudo-orthonormal left-invariant global frame obtained from the pseudo-orthonormal basis of the Lie algebras in Theorem~\ref{th:ce}-(D).

A straightforward calculation shows that gradient vector fields solving the conformally Cotton-flat equation \eqref{eq:hoy3} in case (D.i) are given by $\xi=U_3+\lambda U_4$ for some smooth function satisfying $d\lambda(U_1)=d\lambda(U_2)=d\lambda(U_4)=0$. A direct calculation now shows that, when evaluating on the basis $\{ U_k\}$, one has
$$
\operatorname{Hes}_\varphi=\varphi\left(\begin{array}{cccc}
\lambda &\frac{1}{2}&0&0
\\
\frac{1}{2}&\lambda&0&0
\\
0&0&d\lambda(U_3)+\lambda^2 & \lambda
\\
0&0&\lambda&1
\end{array}\right) \quad\text{and}\quad
\rho= - \left(\begin{array}{cccc}
0&1&0&0
\\
1&0&0&0
\\
0&0&\frac{1}{2}&0
\\
0&0&0&2
\end{array}\right),
$$
from where it follows that $\Delta\varphi=4\lambda\varphi$ and the scalar curvature vanishes. Hence the only non-zero component of the tensor field $\mathfrak{CE}$ is given by
$$
\mathfrak{CE}(U_3,U_3)=\tfrac{1}{2}\varphi(4d\lambda(U_3)+4\lambda^2-1).
$$
This shows that the conformal metric determined by the gradient vector field $\xi=U_3+\lambda U_4$, given by a function $\lambda$ solving  
$$
d\lambda(U_1)=d\lambda(U_2)=d\lambda(U_4)=0
\quad\text{and}\quad
 4d\lambda(U_3)+4\lambda^2=1,
$$ 
is Einstein.

\subsubsection*{Proof of Theorem~\ref{th:ce}-(D.ii)}
Gradient vector fields solving the conformally Cotton-flat equation \eqref{eq:hoy3} in case (D.ii) are given by $\xi=\varepsilon U_3+\lambda U_4$ for some smooth function satisfying $d\lambda(U_1)=d\lambda(U_2)=d\lambda(U_4)=0$. A straightforward calculation now shows that 
$$
\operatorname{Hes}_\varphi\!=\!\varphi\!\left(\!\!\begin{array}{cccc}
\varepsilon(\lambda+1)&\frac{1}{2}\varepsilon\alpha&0&0\\
\frac{1}{2}\varepsilon\alpha&\varepsilon\lambda&0&0\\
0&0&d\lambda(U_3)+\lambda^2&\varepsilon\lambda\\
0&0&\varepsilon\lambda&1
\end{array}\!\!\right),
\,\,
\rho\!=\!-\!\left(\!\!\begin{array}{cccc}
3\varepsilon&\varepsilon\alpha&0&0\\
\varepsilon\alpha&\varepsilon&0&0\\
0&0&\frac{1}{2}\alpha^2+1&\varepsilon\\
0&0&\varepsilon&2
\end{array}\!\!\right),
$$
from where it follows that $\Delta\varphi=\varepsilon(4\lambda+1)\varphi$ and the scalar curvature $\tau=-6\varepsilon$. Now, a straightforward calculation shows that the only non-zero component of $\mathfrak{CE}$ is given by
$$
\mathfrak{CE}(U_3,U_3)=\varphi\left(2d\lambda(U_3)+2\lambda^2-\tfrac{1}{2}\alpha^2-1\right).
$$
This shows that the conformal metric determined by the gradient vector field $\xi=\varepsilon U_3+\lambda U_4$, given by a function $\lambda$ solving 
$$
d\lambda(U_1)=d\lambda(U_2)=d\lambda(U_4)=0
\quad\text{and}\quad
 2d\lambda(U_3)+2\lambda^2-\tfrac{1}{2}\alpha^2=1,
$$ 
is Einstein.

\subsubsection*{Proof of Theorem~\ref{th:ce}-(R)}
Let $\{E_i\}$ be the global orthonormal frame on $H_3\rtimes\mathbb{R}$ obtained  by left-translating the orthonormal basis $\{ e_i\}$ of the Lie algebra. A straightforward calculation shows that the metric is conformal Cotton-flat. Indeed, the vector field $\xi=-(\lambda+\alpha)E_2+\lambda E_4$ solves the equation 
$\operatorname{div}_4W-W(\,\cdot\,,\,\cdot\,,\,\cdot\,,\xi)=0$, and it is a gradient if the smooth function $\lambda$ satisfies $d\lambda(E_1)=d\lambda(E_3)=0$ and $d\lambda(E_2)-d\lambda(E_4)=0$. The Hessian and the Ricci tensors on the global frame $\{ E_i\}$ are expressed as
$$
\begin{array}{rcl}
\operatorname{Hes}_\varphi &=& \varphi\left(\begin{array}{cccc}
\lambda\alpha &0&-\frac{1}{2}\alpha&0\\
0&(\lambda+\alpha)^2-d\lambda(E_4)&0&\lambda(\lambda+\alpha)-d\lambda(E_4)\\
-\frac{1}{2}\alpha &0&\alpha\lambda&0\\
0&\lambda(\lambda+\alpha)-d\lambda(E_4)&0&\lambda^2-d\lambda(E_4)
\end{array}\right)\quad\text{and}
\\  
\noalign{\bigskip}
\rho &= & 
\phantom{\varphi}
\left(\begin{array}{cccc}
2\alpha^2&0&\alpha&0\\
0&-\frac{1}{2}&0&-\frac{1}{2}\\
\alpha&0&2\alpha^2&0\\
0&-\frac{1}{2}&0& -2\alpha^2-\frac{1}{2}
\end{array}\right),
\end{array}
$$
from where it follows that $\Delta\varphi=\alpha(4\lambda+\alpha)\varphi$ and the scalar curvature $\tau=6\alpha^2$.

Hence the non-zero components of $\mathfrak{CE}$ are determined by
$$
\mathfrak{CE}(E_2,E_2)=\mathfrak{CE}(E_4,E_4)=\mathfrak{CE}(E_2,E_4)
= - \tfrac{1}{2}\varphi\left( 4d\lambda(E_4)-4\lambda(\lambda+\alpha) + 1 \right),
$$
which shows that the vector field $\xi=-(\lambda+\alpha)E_2+\lambda E_4$, given by the equations
$$
d\lambda(E_1)=d\lambda(E_3)=0,\quad d\lambda(E_2)-d\lambda(E_4)=0
\quad\text{and}\quad
 d\lambda(E_4)-\lambda(\lambda+\alpha)+\tfrac{1}{4}=0,
$$
is a gradient and the conformal metric induced by the potential function is Einstein.

\subsubsection*{Proof of Theorem~\ref{th:ce}-(L.i)}
Let $\{E_i\}$ be the global orthonormal frame on $H_3\rtimes\mathbb{R}$ obtained  by left-translating the orthonormal basis $\{ e_i\}$ of the Lie algebra.
	The metric is conformally Cotton-flat and,  moreover, gradient vector fields satisfying $\operatorname{div}_4W-W(\,\cdot\,,\,\cdot\,,\,\cdot\,,\xi)=0$  are given by $\xi= \lambda E_3+(\lambda+1)E_4$, 
	for some smooth function $\lambda$ on $H_3\rtimes\mathbb{R}$ 
	such that $d\lambda(E_1)=d\lambda(E_2)=0$ and $d\lambda(E_3)+d\lambda(E_4)=0$.
	A straightforward calculation shows that, on the orthonormal basis $\{ E_i\}$, one has
	$$
	\begin{array}{rcl}
	\operatorname{Hes}_\varphi & =& \varphi\left(\begin{array}{cccc}
	\lambda+1 &\frac{1}{2}\alpha&0&0\\[0.03in]
	\frac{1}{2}\alpha&\lambda+1&0&0\\[0.03in]
	0&0&d\lambda(E_4)+\lambda^2&-d\lambda(E_4)-\lambda(\lambda+1)\\[0.03in]
	0&0&-d\lambda(E_4)-\lambda(\lambda+1)&d\lambda(E_4)+(\lambda+1)^2
	\end{array}\right) 
	\quad\text{and}
	\\
	\noalign{\bigskip}
	\rho&=& -\left(\begin{array}{cccc}
	2&\alpha&0&0\\[0.03in]
	\alpha&2&0&0\\[0.03in]
	0&0&\frac{1}{2}\alpha^2&-\frac{1}{2}\alpha^2\\[0.03in]
	0&0&-\frac{1}{2}\alpha^2&\frac{1}{2} \alpha^2 + 2
	\end{array}\right) ,
	\end{array}
	$$
	from where it follows that $\Delta\varphi=( 4\lambda+3)\varphi$ and the scalar curvature $\tau=-6$. Hence the only non-zero components of the tensor field $\mathfrak{CE}$ are determined by
	$$
	\mathfrak{CE}(E_3,E_3)=\mathfrak{CE}(E_4,E_4)=-\mathfrak{CE}(E_3,E_4)
	= \tfrac{1}{2}\varphi \left( 4d\lambda(E_4) +4\lambda(\lambda+1) - \alpha^2 \right).
	$$
Thus we conclude  that the conformal metric determined by the gradient vector field $\xi= \lambda E_3+(\lambda+1)E_4$, given by a solution of the equations
$$
d\lambda(E_1)=d\lambda(E_2)=0,\quad 
d\lambda(E_3)+d\lambda(E_4)=0
\quad\text{and}\quad
 d\lambda(E_4)+\lambda(\lambda+1)-\tfrac{1}{4}\alpha^2=0,
$$
is Einstein.

\subsubsection*{Proof of Theorem~\ref{th:ce}-(L.ii)}
Let $\{U_i\}$ be the pseudo-orthonormal global frame obtained by left-translating the vectors $\{u_i\}$ on the Lie algebra.
	The metrics in Theorem~\ref{th:ce}-(L.ii) are
	conformally Cotton-flat and, moreover, gradient vector fields satisfying $\operatorname{div}_4W-W(\,\cdot\,,\,\cdot\,,\,\cdot\,,\xi)=0$ are given by  $\xi=\lambda U_2 +4 U_4$, for some smooth function $\lambda$ on $H_3\rtimes\mathbb{R}$ satisfying
	$d\lambda(U_2)=d\lambda(U_3)=0$ and $d\lambda(U_4)=-4\lambda$. The Hessian and the Ricci tensors on the basis $\{ U_i\}$ are given by
	$$
	\operatorname{Hes}_\varphi=\varphi\left(\!\!\begin{array}{cccc}
	d\lambda(U_1) +\lambda^2 + 4\alpha &16&0&0\\
	16&0&0&0\\
	0&0&0&0\\
	0&0&0&16
	\end{array}\!\!\right)
	\quad\text{and}\quad
	\rho= - \left(\!\!\begin{array}{cccc}
	8\alpha& 32&0&0\\
	 32&0&0&0\\
	0&0&0&0\\
	0&0&0& 32
	\end{array}\!\!\right),
	$$
	from where one has that $\Delta\varphi=48\varphi$ and the scalar curvature $\tau=-96$.
	Hence the only non-zero component of $\mathfrak{CE}$ is given by
	$\mathfrak{CE}(U_1,U_1)=2\varphi( d\lambda(U_1) + \lambda^2)$.
	Therefore the conformal metric determined by the gradient vector field $\xi=\lambda U_2 +4 U_4$,
	given by a solution of the equations 
	$$
	d\lambda(U_2)=d\lambda(U_3)=0,\quad d\lambda(U_4)=-4\lambda
	\quad\text{and}\quad
	 d\lambda(U_1)=-\lambda^2,
	$$
	is Einstein.

\subsubsection*{Proof of Theorem~\ref{th:ce}-(L.iii)}
	The Ricci operator of the left-invariant metric in Theorem~\ref{th:ce}-(L.iii) is diagonalizable,  $\Ricci=\operatorname{diag}[-6$, $-6$, $-24$, $-18]$, and therefore $\tau=-54$.
The Weyl curvature operator acting on the space of two-forms has eigenvalues $\{-4,2,2,-4,2,2\}$, where the eigenvalue $2$ is a double root of the minimal polynomial, and hence the metric is weakly generic.
	The metric is also conformally Cotton-flat with left-invariant gradient vector field $\xi= 3 U_4$ satisfying $\operatorname{div}_4W-W(\,\cdot\,,\,\cdot\,,\,\cdot\,,\xi)=0$. Hence it is conformally Einstein by using the results in \cite{KTN}. 
	
	Alternatively one may proceed as in the previous cases
	and consider all gradient solutions of \eqref{eq:hoy3} given by $\xi=\lambda U_4$ for any smooth function satisfying $d\lambda(U_1)=d\lambda(U_2)=d\lambda(U_3)=0$. Hence one has
	$$
	\operatorname{Hes}_\varphi=\varphi\left(\begin{array}{cccc}
	0&\lambda&0&0\\\lambda&0&0&0\\0&0&4\lambda&0\\0&0&0&d\lambda(U_4)+\lambda^2
	\end{array}
	\right)
	\quad\text{and}\quad
	\rho= - \left(\begin{array}{cccc}
	0&6&0&0\\
	6&0&0&0\\
	0&0&24&0\\
	0&0&0&18
	\end{array}\right)	,
	$$
	so that $\Delta\varphi= \left( d\lambda(U_4)+\lambda(\lambda+6) \right) \varphi$. 
	Hence the only non-zero terms in $\mathfrak{CE}$ are determined by
	$$
	\begin{array}{l}
	\mathfrak{CE}(U_1,U_2)=-\frac{1}{2}\varphi\left( d\lambda(U_4)+(\lambda+5)(\lambda-3)\right),
	\\
	\noalign{\medskip}
	\mathfrak{CE}(U_3,U_3)=-\frac{1}{2}\varphi\left( d\lambda(U_4)+(\lambda-3)(\lambda-7)\right),
	\\
	\noalign{\medskip}
	\mathfrak{CE}(U_4,U_4)=\phantom{-}\frac{3}{2}\varphi\left( d\lambda(U_4)+(\lambda+1)(\lambda-3)\right),	
	\end{array}
	$$
	from where it follows that the function $\lambda$ is constant $\lambda=3$.

\begin{remark}\rm
A straightforward calculation shows that all the conformally Einstein metrics in Theorem~\ref{th:ce} are indeed conformally Ricci-flat. While the conformally Einstein metric in  Theorem~\ref{th:ce}-(L.iii) is  unique (up to scaling), the other cases admit infinitely many conformally Einstein metrics, from where it follows that they are conformally equivalent to a Ricci-flat $pp$-wave (see \cite{Bri}).
\end{remark}

\section{Proof of Theorem~\ref{th:strictly Bach flat}}\label{se:proof-2}
Considering the left-invariant metric in Theorem~\ref{th:strictly Bach flat}-(i), a straightforward calculation shows that the Weyl curvature operator $W:\Lambda^2\rightarrow\Lambda^2$ has six-distinct non-zero complex eigenvalues and thus it is weakly generic. Next, we show that the metric does not satisfies the conformal $C$-space condition, i.e., there does not exist a (not necessarily gradient) vector field $X$ so that 
$\operatorname{div}_4W-W(\,\cdot\,,\,\cdot\,,\,\cdot\,,X)=0$. Hence the metric  is strictly Bach-flat.

Let $\{ V_i\}$ be the orthonormal frame on $H_3\rtimes\mathbb{R}$ obtained by left-translating the vectors $\{ v_i\}$ at the Lie algebra  and set $X=\sum X^j V_j$. Now a straightforward calculation shows that
$$
\begin{array}{l}
\operatorname{div}_4W(V_4,V_1,V_1)-W(V_4,V_1,V_1,X)=26X^4,
\\
\operatorname{div}_4W(V_4,V_1,V_2)-W(V_4,V_1,V_2,X)=4\sqrt{14} X^3,
\\
\operatorname{div}_4W(V_4,V_1,V_3)-W(V_4,V_1,V_3,X)=-4\sqrt{14} X^2,
\\
\operatorname{div}_4W(V_4,V_1,V_4)-W(V_4,V_1,V_4,X)=54\sqrt{14}-26 X^1,
\\
\operatorname{div}_4W(V_3,V_2,V_4)-W(V_3,V_2,V_4,X)=144-4\sqrt{14}X^1,
\end{array}
$$
from where one has that the above equations have no solution, and thus $H_3\rtimes\mathbb{R}$ is not a conformal $C$-space.

We proceed in an analogous way with metrics in Theorem~\ref{th:strictly Bach flat}-(ii). The Weyl curvature operator acting on the space of two-forms is three-step nilpotent and it is not a  conformal $C$-space, which shows that these metrics are also strictly Bach-flat. We omit the details that are completely analogous to the previous case.

\section{Conclusions}

Bach-flatness is a very restrictive condition for left-invariant Riemannian metrics on four-dimensional Lie groups \cite{AGS,CL-GM-GR-GR-VL}. 
The Lorentzian situation is more subtle due to the fact that the restriction of the metric to the three-dimensional normal subgroup $G$ may be a positive definite, Lorentzian or degenerate metric. 

We classify all left-invariant Bach-flat Lorentzian metrics on semi-direct extensions $H_3\rtimes\mathbb{R}$. As a consequence it is shown that the class of conformally Einstein metrics which are not locally conformally flat reduces to plane waves and six generically non-homothetic classes (see Theorem~\ref{th:ce}). 
On the opposite, the class of strictly Bach-flat metrics, i.e., those which are not conformally Einstein, reduces to two non-homothetic families (see Theorem~\ref{th:strictly Bach flat}). 

A special situation occurs when the semi-direct extension is a product.
Bach-flat Lorentzian metrics on the product Lie group $H_3\times\mathbb{R}$ are locally conformally flat or a plane wave with parallel Ricci tensor (see Remark~\ref{re:producto}).

Among the Bach-flat metrics there is a single one which is critical for all quadratic curvature functionals without being neither symmetric nor a plane wave (see Theorem~\ref{th:ce}--(D.i)). This is in sharp contrast with the Riemannian situation \cite{criticas}. Moreover, the conformally Einstein metric in Theorem~\ref{th:ce}--(L.iii) is a shrinking algebraic  soliton for the RG2 flow, while the family of metrics in Theorem~\ref{th:ce}--(L.ii) provides Bach-flat steady solitons for the RG2 flow, which are therefore self-similar solutions of the flow.

Finally, note that all conformally Einstein semi-direct extensions $H_3\rtimes\mathbb{R}$ are conformal to a $pp$-wave except metrics corresponding to Theorem~\ref{th:ce}--(L.iii) whose Weyl curvature operator acting on the space of two-forms is not nilpotent.


\begin{thebibliography}{99}
\bibitem{AGS}
E. Abbena, S. Garbiero, and S. Salamon,
Bach-flat Lie groups in dimension $4$, 
\emph{C. R. Math. Acad. Sci. Paris} \textbf{351} (2013), 303--306.

\bibitem{Bach}
R. Bach, 
Zur Weylschen Relativit\"atstheorie und der Weylschen Erweiterung des Kr\"ummungstensorbegriffs,
\emph{Math. Z.} \textbf{9} (1921), 110--135.

\bibitem{BO03} 
M. Blau and M. O'{}Loughlin, Homogeneous plane waves
\emph{Nuclear Physics B}  \textbf{654} (2003), 135--176.

\bibitem{BPR} 
H. Bondi, F. A. E. Pirani, and I. Robinson, 
Gravitational waves in general Relativity III. Exact plane waves, 
\emph{Proc. R. Soc. Lond. Ser. A Math. Phys. Eng. Sci.} \textbf{251} (1959), 519--533.

\bibitem{Bri1}
H. W Brinkmann, 
Riemann spaces conformal to Einstein spaces, 
\emph{Math. Ann.} \textbf{91} (1924), 269--278.

\bibitem{Bri}
H. W Brinkmann, 
Einstein spaces which are mapped conformally on each other, 
\emph{Math. Ann.} \textbf{94} (1925), 119--145.

\bibitem{criticas}
M. Brozos-Vázquez, S. Caeiro-Oliveira, E. García-Río, and R. Vázquez-Lorenzo,
Four-dimensional homogeneous critical metrics for quadratic curvature functionals,
to appear.

\bibitem{BGV19} 
M. Brozos-V\'azquez, E. Garc\'ia-R\'io, and X. Valle-Regueiro, Isotropic quasi-Einstein manifolds,
\emph{Classical Quantum Gravity}  \textbf{36} (2019), 245005 (13pp).


\bibitem{Cahen}
M. Cahen, J. Leroy, M. Parker, F. Tricerri, and L. Vanhecke,
Lorentz manifolds modelled on a Lorentz symmetric space,
\emph{J. Geom. Phys.} \textbf{7} (1990), 571--581. 

\bibitem{CaZa2}
G. Calvaruso and A. Zaeim, 
Conformal geometry of semi-direct extensions of the Heisenberg group,
\emph{J. Math. Phys. Anal. Geom.} \textbf{17} (2021), 407--421.

\bibitem{CaZa}
G. Calvaruso and A. Zaeim, 
Four-dimensional Lorentzian Lie groups,
\emph{Differential Geom. Appl. } {\bf 31} (2013), 496--509. 

\bibitem{CC}
G. Calvaruso and M. Castrillón,
Cyclic Lorentzian Lie groups, 
\emph{Geom. Dedicata} \textbf{181} (2016), 119--136.

\bibitem{CL-GM-GR-GR-VL}
E. Calviño-Louzao, X. Garc\'ia-Mart\'inez, E. Garc\'ia-R\'io, I. Guti\'errez-Rodr\'iguez, and R. V\'azquez-Lorenzo, 
Conformally Einstein and Bach-flat four-dimensional homogeneous manifolds,
\emph{J. Math. Pures Appl. (9)} \textbf{130} (2019), 347--374.

\bibitem{Carfora}
M. Carfora,
Renormalization group and the Ricci flow, 
\emph{Milan J. Math.} \textbf{78} (2010), 319--353.

\bibitem{CoHePe}
A. Coley, S. Hervik, and N. Pelavas,
Spacetimes characterized by their scalar curvature invariants,
\emph{Classical Quantum Gravity} \textbf{26} (2009), 025013 (33 pp).

\bibitem{Cox} 
D. Cox, D. Little, and D. O'Shea,
\emph{Ideals, varieties, and algorithms. An introduction to computational algebraic geometry and commutative algebra},
Undergraduate Texts in Mathematics. Springer, Cham, 2015.

\bibitem{Singular}
W. Decker, G.-M. Greuel, G. Pfister, and H. Sch\"onemann,
\emph{{\sc Singular} {4-3-0} --- {A} computer algebra system for polynomial computations}, \texttt{https://www.singular.uni-kl.de}, 2022.

\bibitem{13}
J. Ehlers and W. Kundt, Exact solutions of the gravitational field equations, 
\emph{Gravitation: an introduction to current research}, 49--101, Wiley, New York, 1962.

\bibitem{GGI2}
K. Gimre, Ch. Guenther, and J. Isenberg,
A geometric introduction to the two-loop renormalization group flow,
\emph{J. Fixed Point Theory Appl.} \textbf{14} (2013), 3--20.

\bibitem{GN}
A. R. Gover and P. A. Nagy, 
Four-dimensional conformal $C$-spaces,
\emph{Q. J. Math.} \textbf{58} (2007), 443--462.

\bibitem{KT}
Y. Kondo and H. Tamaru,
A classification of left-invariant Lorentzian metrics on some nilpotent Lie groups, 
\emph{Tohoku Math. J. (2)} \textbf{75} (2023), to appear.

\bibitem{KTN}
C. N. Kozameh, E. T. Newman, and K. P. Tod,
Conformal Einstein spaces, 
\emph{ Gen. Relativity Gravitation} \textbf{17} (1985), 343--352.
	
\bibitem{Kulkarni}
R. S. Kulkarni,
Curvature and metric,
\emph{Ann. of Math. (2)} \textbf{91} (1970), 311--331.

\bibitem{Lauret}
J. Lauret,
Ricci soliton homogeneous nilmanifolds, 
\emph{Math. Ann.} \textbf{319} (2001), 715--733.

\bibitem{Leistner}
T. Leistner,
Conformal holonomy of C-spaces, Ricci-flat, and Lorentzian manifolds,
\emph{Differential Geom. Appl.} \textbf{24} (2006), 458--478.


\bibitem{LN10} 
	T. Leistner and P. Nurowski, Ambient Metrics for $n$-dimensional $pp$-waves
	\emph{Comm. Math. Phys.}  \textbf{296} (2010), 881--898. 

\bibitem{Mannheim1}
P. D. Mannheim, 
Alternatives to dark matter and dark energy,
\emph{Prog. Part. Nucl. Phys.} \textbf{56} (2006), 340--445. 


\bibitem{Mannheim2}
P. D. Mannheim, 
Making the case for conformal gravity,
\emph{Found. Phys.} \textbf{42} (2012), 388--420.

\bibitem{Milnor}
J. Milnor,
Curvatures of left invariant metrics on Lie groups,
\emph{Adv. Math.} \textbf{21} (1976), 293--329.

\bibitem{MuRi}
D. Müller and F. Ricci,
Analysis of second order differential operators on Heisenberg groups, I,
\emph{Invent. math.} \textbf{101} (1990), 5454--582.

\bibitem{Rahmani}
S. Rahmani,
Métriques de Lorentz sur les groupes de Lie unimodulaires, de dimension trois,
\emph{J. Geom. Phys.} \textbf{9} (1992), 295--302. 

\bibitem{St}
R. F. Streater,
The representations of the oscillator group,
\emph{Comm. Math. Phys.} \textbf{4} (1967), 217--236.

\bibitem{Wears}
Th. H. Wears,
On algebraic solitons for geometric evolution equations on three-dimensional Lie groups,
\emph{Tbilisi Math. J.} \textbf{9} (2016), 33--58.

\end{thebibliography}
\end{document}